\documentstyle[12pt,amssymb,mysection]{article}
\voffset = - 2.3 truecm
\newcommand{\boxx}{\ \vrule width 2.0 mm height 2.0 mm  depth 0 mm }

\newcounter{theo}[section]
\newcounter{lemma}[section]
\newcounter{defi}[section]
\newcounter{rem}[section]
\newcounter{cor}[section]
\catcode`@=11
\catcode`@=12
\def\thetheo{\thesection.\arabic{theo}}
\def\thelemma{\thesection.\arabic{lemma}}
\def\thedefi{\thesection.\arabic{defi}}
\def\therem{\thesection.\arabic{rem}}
\def\thecor{\thesection.\arabic{cor}}
\newenvironment{theo}{%
\refstepcounter{theo} {THEOREM \thetheo.}\bgroup\rm}{\egroup}
\newenvironment{lemma}{%
\refstepcounter{lemma} {LEMMA \thelemma.}\bgroup\rm}{\egroup}
\newenvironment{defi}{%
\refstepcounter{defi} {DEFINITION \thedefi.}\bgroup\rm}{\egroup}
\newenvironment{rem}{%
\refstepcounter{rem} {REMARK \therem.} \bgroup\rm}{\egroup}

\newfont{\cmss}{cmss12 scaled 1000}
\newfont{\cmex}{cmex10 scaled 1440}
\newfont{\eufm}{eufm10 scaled 1200}
\newcommand{\Int}{\mbox{\cmex\symbol{90}}}
\newcommand{\Rgot}{\mbox{\eufm\symbol{82}}}
\newcommand{\Cgot}{\mbox{\eufm\symbol{67}}}
\begin{document}
\ \vskip 3cm
\begin{minipage}{14.0cm}
{\bf ON TRANSFORMATIONS OF POTAPOV's \\
FUNDAMENTAL MATRIX INEQUALITY \\[0.5cm]
V.E. Katsnelson  
}
\end{minipage}
\vskip .5cm

\begin{abstract}
According to V.P.Potapov, a classical interpolation problem can be
 reformulated in terms  of a so-called Fundamental Matrix Inequality
 (FMI).  To show that every solution of the FMI satisfies the interpolation
 problem, we usualy have to transform the FMI in some
special way.
 In this paper the number of of  transformations of the FMI
which come into play are motivated and
 demonstrated by simple,
but typical examples.
\end{abstract}
\setcounter{section}{-1}
\begin{minipage}{15.0cm}
\section{\hspace{-0.4cm}.\hspace{0.2cm} PREFACE}
\end{minipage}\\[-0.2cm]
\setcounter{equation}{0}

V.P.Potapov's approach to classical interpolation problems research consists
in the following. Instead of original interpolation problem (or  problem on
integral representation), an inequality for analytic functions 
is considered in an appropriate domain. This inequality is said to be 
{\sl the Fundamental Matrix Inequality} (FMI) for the considered interpolation
 problem. Here two problems appear. The first problem is  how to ``solve'' this
inequality. The second problem is to prove that this inequality is equivalent
to the original interpolation problem. 

The study of the second problem consist of two parts. First,
 we have to prove that any function which is a solution of the original
 problem is also a solution of the  FMI. Usually this part is not difficult.
Secondly, we have to extract the full interpolation information from the 
FMI. This means that we have to prove that any analytic function which
 satisfies the FMI  is also a solution of of the original interpolation
 problem.
In simple situations it is not difficult to obtain the interpolation
 information from the FMI. However, in the general case this is not easy,
and {\ we have to apply a special transformation to the FMI}.
 Such a transformation can be applied to  every FMI. However, in the
 simplest situations it is possible to do without such a transformation.
The development of Potapov's method began with consideration of the
 simplest interpolation problem, i.e. the Nevanlinna-Pick ($\cal NP$) problem. 
The equivalence of the $\cal NP$ problem to its FMI   is clear.
Because of this, this transform was camouflaged in the 
beginning of the theory. However, by the study of of the power
moment problem we already can not do without it. In the paper \cite{KKY} such a
transform was used in the very general setting of the so called {\sl
Abstract Interpolation Problem}. Namely, such a transformation was used
in considerations related to Theorem 1 of this paper. Of course, the
authors of \cite{KKY} took into account the experience which was accumulated
 by  previous work with concrete problems. However, this transformation
was introduced in \cite{KKY} in a formal way, without any motivation.
As result, the proof of Theorem 1 of \cite{KKY} looks like a trick.
This is not satisfactory, because the transformation of FMI  lies at the heart
of the FMI  business. The main goal of the present paper is to motivate and
to demonstrate the transformation of the FMI by the simplest but typical
 example of the power moment problem. For  contrast, the $\cal NP$
 problem and the FMI  for it are
 considered as well. We would like to demonstrate the algebraic side of
the matter. Therefore, we will avoid the entourage of general vector spaces and
 Hilbert spaces in the generality of the paper \cite{KKY}. All our spaces are
finite-dimensional. Instead of abstract kernels and operators, we will consider
matrices. 

\vspace*{0.2cm} 
\begin{minipage}{15.0cm}
\section{\hspace{-0.4cm}.\hspace{0.2cm}THE FMI  AND ITS STRUCTURE}
\end{minipage}\\[-0.2cm]
\setcounter{equation}{0}

Classical interpolation problems can be considered for various function
classes in various domains. Here
we consider two function classes related to the unit disc $\Bbb D$ and to
the upper half plane $\Bbb H$.

\begin{defi}\ \ 
I.{\sl  {\cmss The class} {\mbox{\Cgot}\,($\Bbb D$)} is the class
of functions $w$ which are holomorphic outside the unit circle 
$\Bbb T$, satisfy the symmetry condition
\begin{equation}
 w(z) =- w^{*}(1/\overline{z}) \qquad (\,  z \in {\Bbb C}\setminus{\Bbb T}\, ) 
\label{1.1}
\end{equation}         
and the positivity condition
\begin{equation}
\frac{w(z)+w^{*}(z)}{1-|z|^{2}}\geq 0
\qquad (\,  z \in {\Bbb C}\setminus{\Bbb T} \, ) .
\label{1.2}
\end{equation}
}

II.{\sl {\cmss The class {\mbox{\Rgot}\,($\Bbb H$)}} is the class of
 functions $w$ which
 are holomorphic outside the real axes 
${\Bbb R}$ and satisfies the symmetry condition
\begin{equation}
w(z)=w^{*}(\overline{z})
\quad (\,  z \in {\Bbb C}\setminus{\Bbb R} \, ) 
\label{1.3}
\end{equation}
and the positivity condition
\begin{equation}
\frac{w(z)-w^{*}(z)}{z-\overline{z}}\geq 0
\qquad (\,  z \in {\Bbb C}\setminus{\Bbb R} \, ).
\label{1.4} 
\end{equation}
}
III. {\sl
 {\cmss
 The class {$\mbox{\Rgot}_0$\,($\Bbb H$)}
}
 is the subclass of the class {\mbox{\Rgot}\,($\Bbb H$)} which is
singled out by the condition
\begin{equation}
\overline{\lim_{y \, \uparrow\, \infty }}\,y\,|w(iy)| <\infty .
\label{1.4.a}
\end{equation}
}
\end{defi}

The FMI of a classical interpolation problem has the form
\begin{equation}
\left[
\begin{array}{ccc}
A& \hspace*{-0.3cm} & B_{w}(z)\cr
 \vspace*{-0.1cm}&\vspace*{-0.1cm}\hspace*{-0.3cm} &\vspace*{-0.1cm} \cr
B^{\raisebox{0.3ex}{\mbox{$\ast$}}}_{w}(z)&\hspace*{-0.3cm} &C_{w}(z)
\end{array}
\right]
\geq 0,
\label{1.5}
\end{equation}
where $A$ is some hermitian matrix, constructed from the interpolation data
(interpolation points and interpolating values) only. It is nonnegative if
and only if the considered interpolation problem is solvable.
The entry $C_{w}(z)$ contains the function $w$ only, but not the
interpolation data. Its form depend on the function class to which the
 function $w$ belongs. For an interpolation problem in the class 
\Cgot\,($\Bbb D$)
 the entry $C_{w}(z)$ has the form
\begin{equation}
C_{w}(z)=\frac{w(z)+w^{\ast}(z)}{1-|z|^{2}}.
\label{1.6}
\end{equation}
 For an interpolation problem in the class \Rgot\,($\Bbb H$)
 the entry $C_{w}(z)$ has the form
\begin{equation}
C_{w}(z)=\frac{w(z)-w^{\ast}(z)}{z-\bar{z}}.
\label{1.7}
\end{equation}
In the entry $B_{w}(z)$ both the interpolation data and the function $w$ are
combined. This entry looks like
\begin{equation}
B_{w}(z)=(zI-T)^{-1}(u\cdot w(z) - v),
\label{1.8}
\end{equation}
or like   
\begin{equation}
B_{w}(z)=T(I-zT)^{-1}(u\cdot w(z) - v)
\label{1.9}
\end{equation}
To each classical interpolation problem the following objects are
 related:\\[0.2cm]
\hspace*{1.0cm} \begin{minipage}{14.0cm}
1.The hermitian matrix $A$, which is nonnegative iff the problem is solvable.\\
2. The matrix $T$ which ``determines'' the interpolation nodes.\\
3.The vectors $u$ and $v$ which determine the interpolation values.\\
\end{minipage}

The terms $A, T, u, v$ satisfy the so called {\sl Fundamental Identity} ( FI).
 The
form of the FI depends on the function class in which the interpolation 
problem is considered. For the function class \Cgot\,($\Bbb D$), FI
has the form 
\begin{equation}
A - TAT^{\ast}=uv^{\ast}+vu^{\ast}.
\label{1.10}
\end{equation}
For the class \Rgot\,($\Bbb H$), FI has the form 
\begin{equation}
TA - AT^{\ast}=uv^{\ast} -vu^{\ast}.
\label{1.11}
\end{equation}
 If the FMI (\ref{1.5}) is satisfied (for some $z$), and if $M$ is a matrix
of an appropriate size, then the inequality
\begin{equation}
  M\ 
\left[
\begin{array}{ccc}
A& \hspace*{-0.3cm} & B_{w}(z)\cr
 \vspace*{-0.1cm}&\vspace*{-0.1cm}\hspace*{-0.3cm} &\vspace*{-0.1cm} \cr
B^{\raisebox{0.3ex}{\mbox{$\ast$}}}_{w}(z)&\hspace*{-0.3cm} &C_{w}(z)
\end{array}
\right]
\ M^{\ast}\geq 0
\label{1.12}             
\end{equation}
holds as well. If the matrix $M$ is invertible, then both the inequalities
(\ref{1.5}) and (\ref{1.12})  are equivalent.

\vspace*{0.2cm} 
\begin{minipage}{15.0cm}
\section{\hspace{-0.4cm}.\hspace{0.2cm} FMI FOR THE NEVANLINNA --
 PICK PROBLEM.}
\end{minipage}\\[-0.2cm]
\setcounter{equation}{0}

Now we obtain the FMI  for the $\cal NP$ problem in the function class
\Cgot\,($\Bbb D$).

\begin{defi}
{\sl Given $n$ points $z_1,\,z_2,\, \ldots,\,z_n$
in the unit disc $\Bbb D$  {\cmss (interpolation nodes)}   and $n$ complex
numbers $w_1,\,w_2,\, \ldots,\,w_n$  {\cmss (interpolation values)}. A
 holomorphic function $w(z)$ from the class \Cgot\,($\Bbb D$)  is said to
be {\cmss  a solution of the Nevanlinna -- Pick problem with  interpolation
 data
$\{z_1,w_1\}\, , \{z_2,w_2\}\, , \ldots \, ,\{z_n,w_n\} $}, if the
{\cmss  interpolation conditions} 
\begin{equation}
      w(z_k) = w_k \qquad ( k=1, 2, \ldots , n )
\label{2.1}
\end{equation}
are satisfied.}
\end{defi}

Let us associate with the $\cal NP$  problem two $n\times 1$ vectors, which
 characterize
the interpolation values:
\begin{equation}
u = 
\left[
\begin{array}{c}
1 \cr
\vdots \cr
1 \cr
\end{array}
\right] \qquad \mbox{and} \qquad
v = 
\left[
\begin{array}{c}
w_1 \cr
\vdots \cr
w_n \cr
\end{array}
\right].
\label{2.2}
\end{equation}
The matrix $T$, which characterize the interpolation nodes, has the form
\begin{equation}
T= {\rm diag}\: [\, z_1,\, z_2\,, \cdots\, ,  z_n\, ]\, .
\label{2.3}
\end{equation}
The matrix $A$, the so called {\sl Pick matrix}  for the problem,  has
 the form
\begin{equation}
A=\big\| a_{kl}\big\|_{1\leq k,l \leq n}\, ,\qquad 
a_{kl}= \frac{w_k+\bar{w}_l}{1-z_k\bar{z}_l}.
\label{2.4}
\end{equation}
The Fundamental Identity (\ref{1.10}) for this chois of
 $u,\, v,\, T, \mbox{and} A $ can be
checked directly.\\

The Fundamental Matrix Inequality for the Nevanlinna-Pick problem
( FMI({$\cal NP$}) ) has the form (\ref{1.5}) with $A$ from (\ref{2.4}), 
$C_{w}(z)$ from (\ref{1.6}) and  $B_{w}(z)$ from
 (\ref{1.8}),\,(\ref{2.2})\,(\ref{2.3}). 

\begin{theo} ({\cmss
From {\rm FMI($\cal NP$)} to interpolation conditions.})
{\sl Let $w(z)$ be a function which is 
 holomorphic in the unit disc $\Bbb D$  and which satisfies the
 FMI({$\cal NP$})  for every $z \in \Bbb D$.
Then the function $w$ satisfies the condition
$w(z)+w^{\ast}(z)\geq 0\; (z\in \Bbb D)$ and
the interpolation conditions {\rm (\ref{2.1})}.}
\end{theo}

PROOF. Since the entry $C_{w}(z)$ must be nonnegative for $z\in \Bbb D$, the
 real part of the function $w$ is nonnegative
in\footnote{If we continue the function $w$,
which is defined originally in $\Bbb D$ only,
into the exterior of the unit circle according to the symmetry (\ref{1.1}),
then the function which is  continued in this way  will satisfy
 the condition (\ref{1.2}).}
 $\Bbb D$.
Now we take into account the concrete
form of the entry $B_{w}(z)$:
\begin{equation}
B_w(z)=
\left[
\begin{array}{c}
b_{1,w}(z)\cr
b_{2,w}(z)\cr
\vdots\cr
b_{n,w}(z)\cr
\end{array}
\right],
\label{2.5}
\end{equation}
where
\begin{equation}
b_{k,w}(z)= \frac{w(z)-w_k}{z-z_k}\quad (k=1,2,\cdots,n).
\label{2.6}
\end{equation}
Because the ``full'' matrix (\ref{1.5}) is nonnegative, its appropriate 
submatrices are nonnegative all the more:
\begin{equation}
\left[
\begin{array}{ccc}
a_{kk}& \hspace*{-0.3cm} & b_{w, k}(z)\cr
 \vspace*{-0.1cm}&\vspace*{-0.1cm}\hspace*{-0.3cm} &\vspace*{-0.1cm} \cr
b^{\raisebox{0.3ex}{\mbox{$\ast$}}}_{w,k}(z)&\hspace*{-0.3cm} &C_{w}(z)
\end{array}
\right]
\geq 0,
\label{2.7}
\end{equation}
 Since the function $w$ is holomorphic in $\Bbb D$, the entry
$C_{w}(z)$ , (\ref{1.6}) , is locally bounded in $\Bbb D$.
Thus, from (\ref{2.7}) it follows, that the entry $b_{w,k}(z)$
is locally bounded in $\Bbb D$ as well. However,if  function $b_{k}$ is
 bounded even near the point $z_{k}$, then
the interpolation conditions (\ref{2.1}) are satisfied. 
\hfill $\boxx$\\[0.2cm]
Thus, for the $\cal NP$ interpolation problem it is not difficult to extract
 the
 interpolation information from its FMI.

It is worth mentioning, that the inequality (\ref{2.7}) can be consider as
an inequality of the form (\ref{1.12}), with
\begin{equation}
\raisebox{-0.3cm}{$M=$}
\begin{array}{c}
\mbox{\scriptsize \hspace{2.3cm}k\hspace{1.8cm} n+1\hspace{0.1cm}} \cr
\left[
\begin{array}{cccccccc}
0&0 &\cdots &1\cdots&0&\vdots& 0\cr
0&0 &\cdots &0\cdots&0&\vdots& 1\cr
\end{array}
\right]\cr
\end{array}\: \: \raisebox{-0.3cm}{$\cdot$}
\label{2.8}
\end{equation}

\vspace*{0.2cm} 
\begin{minipage}{15.0cm}
\section{\hspace{-0.4cm}.\hspace{0.2cm}  DERIVATION  OF THE
 FMI\,\boldmath$( \cal NP)$\unboldmath }
\end{minipage}\\[-0.2cm]
\setcounter{equation}{0}

A crucial role in deriving of the FMI for the
$\cal  NP$ problem is played by the Riesz-Herglotz 
theorem. Given a nonnegative measure $\sigma$ and a real number $c$,
we associate with them the function $w_{\sigma, c}$:
\begin{equation}
w_{\sigma, c}(z)=ic+ \frac{1}{2}
\int\limits_{\Bbb T}\frac{t+z}{t-z}\,d\sigma (t) ,
 \qquad ( z \in \Bbb C \setminus \Bbb T ).
\label{3.1}
\end{equation}
The function $w_{\sigma, c}$ belongs to the class \Cgot\,($\Bbb D$).

   THEOREM ({\cmss RIESZ-HERGLOTZ}). {\sl  Let $w$ be a function which
belongs to the class \Cgot\,($\Bbb D$). Then this function $w$ is of the form
{\rm (\ref{3.1})} for some $\sigma$ and $c$. Such $\sigma$ and $c$ are
determined from the given 
$w$ uniquely.}

Let us start to derive the FMI($\cal NP$). Given a measure $\sigma \geq 0$ on 
$\Bbb T$, a real number $c$ and points $z_1, z_2, \cdots, z_n; z \in \Bbb D,$.
Let $u$ be defined by (\ref{2.2}), $T$ be defined by (\ref{2.3}). Then
the following inequality ($z_1, z_2, \cdots, z_n$ appear in $T$) holds:
\begin{equation}
\raisebox{0.65cm}{$\Int$ \hspace{-0.55cm}
\raisebox{-1.4cm}{$\scriptstyle\Bbb T$}}\hspace{0.25cm}
\left[
\begin{array}{c}
\vspace*{-0.1cm}
(tI-T)^{-1}u \cr
\vspace*{-0.15cm}
\cdot - \cdot - \cdot - \cdot 
\vspace*{0.15cm}\cr 
\displaystyle \bar{t}(\bar{t}-\bar{z})^{-1}\cr 
\end{array}
\right]\,
\cdot
d\sigma(t)
\cdot
\left[
\begin{array}{ccc}
u^{\ast}(\bar{t}I-T^{\star})^{-1}&\hspace{-0.2cm} 
\begin{array}{c}
|\vspace*{-0.1cm}\cr
\vspace*{-0.1cm}
\cdot
 \cr
\vspace*{-0.1cm}|\vspace*{0.1cm}\cr
\end{array}
\hspace{-0.2cm}  &\displaystyle \frac{t}{t-z}
\end{array}
\right]\,
\geq 0 . 
\label{3.2}
\end{equation}
This is a block-matrix inequality of the form
\begin{equation}
 \left[
\begin{array}{ccc}
A_{\sigma}& \hspace*{-0.3cm} & B_{\sigma}(z)\cr
 \vspace*{-0.1cm}&\vspace*{-0.1cm}\hspace*{-0.3cm} &\vspace*{-0.1cm} \cr
B^{\raisebox{0.3ex}{\mbox{$\ast$}}}_{\sigma}(z)&\hspace*{-0.3cm} &C_{\sigma}(z)
\end{array}
\right]
\geq 0.
\label{3.3}
\end{equation}
We consider also the function $w_{\sigma,c}$ , (\ref{3.1}) , associated with 
$\sigma$ and $c$.

Now we will discuss the entries of the block-matrix on the right-hand side
of the inequality (\ref{3.3}). Originally these entries were defined by means
of an integral representation. However, they can be expressed in terms of
the function $w_{\sigma, c}$. Let us consider the block $A_{\sigma}$:
\begin{equation}
A_{\sigma}=\int\limits_{\Bbb T}(tI-T)^{-1}u
\cdot
d\sigma (t)
\cdot
u^{\ast}(\bar{t}-T^{\ast})^{-1},
\label{3.4}
\end{equation}
or, for the entries $A_{\sigma}=\|a_{\sigma,kl}\|_{1\leq k,l \leq n}$ :
$$
a_{\sigma,kl}=\int\limits_{\Bbb T}(t-z_{k})^{-1}
\cdot
d\sigma (t)
\cdot
(\bar{t}-\bar{z}_{l})^{-1}, \qquad ( 1\leq k, l \leq n).
$$
According to the well known identity for the {\sl Schwarz kernel}\ \
$\displaystyle\frac{1}{2}\, (t+z)(t-z)^{-1}$,
\begin{equation}
A_{\sigma}=\left\|
\frac
{w_{\sigma , c}(z_k)+\overline{w_{\sigma , c}(z_l)}}{1-z_k\bar{z}_l}
\right\|_{1\leq k,l \leq n}\cdot
\label{3.5}
\end{equation}
(The constant $c$ does not appear in (\ref{3.5}).) The block $B_{\sigma}$
has the following form:
\begin{equation}
B_{\sigma }=\int\limits_{\Bbb T}
(tI-T)^{-1}u\cdot\frac {t}{t- z}\cdot d\sigma(t).
\label{3.6}
\end{equation}
The block $B_{\sigma }$ (which  does not depend on c)
 can be transformed in the following way.
 Integrating the identity
$$
(tI-T)^{-1}\frac{t}{t-z}=(zI-T)^{-1}\cdot \frac{1}{2}\,\frac{t+z}{t-z}-
(zI-T)^{-1}\cdot \frac{1}{2}\frac{tI+T}{tI-T}
$$
with respect the measure $d\sigma$, we obtain:
\begin{equation}
B_{\sigma }=(zI-T)^{-1}\Big(
uw_{\sigma,c}(z)-v_{\sigma,c}
\Big),
\label{3.7}
\end{equation}
where
\begin{equation}
v_{\sigma,c}=icu+\frac{1}{2}
\int\limits_{\Bbb T}\frac{tI+T}{tI-T}\,d\sigma(t).
\label{3.8}
\end{equation}
It can be checked that
\begin{equation}
A_{\sigma}-TA_{\sigma}T^{\ast}
=u\cdot v^{\ast}_{\sigma,c}-v\cdot u^{\ast}_{\sigma,c} .
\label{3.9}
\end{equation}
According to (\ref{3.1}) and to (\ref{2.3}),
\begin{equation}
v_{\sigma,c}=
\left[
\begin{array}{c}
w_{\sigma,c}(z_1)\cr
w_{\sigma,c}(z_2)\cr
\cdots \cr
w_{\sigma,c}(z_n)\cr
\end{array}
\right]
\cdot
\label{3.10}
\end{equation}
Of course,
\begin{equation}
C_{\sigma}(z)=\int\limits_{\Bbb T}
\frac{d\sigma(t)}{|t-z|^2}=
\frac{w_{\sigma,c}(z)+\overline{w_{\sigma,c}(z)}}{1-|z|^2}\,\cdot
\label{3.11}
\end{equation}
Now  let the function $w_{\sigma ,c}$ satisfy the interpolation conditions
 (\ref{2.1}) , i.e. let
\begin{equation}
w_{\sigma, c}(z_k)=w_k     \qquad  (k=1, 2,\, \dots\, , n)\, .
\label{3.12}
\end{equation}
Comparing (\ref{3.5}) and (\ref{2.4}), we obtain that 
\begin{equation}
A_{\sigma}=A.
\label{3.13}
\end{equation}
From (\ref{3.10}) and (\ref{2.2}),
\begin{equation}
 v_{\sigma,c}=v.
\label{3.14}
\end{equation}
 Comparing now (\ref{3.7}) with (\ref{1.8}),  we obtain that
\begin{equation}
B_{\sigma}(z)=B_{w_{\sigma ,c}}(z)\, .
\label{3.15}
\end{equation}
Of course, (\ref{3.11})), $C_{\sigma}(z)=C_{w_{\sigma,c}}(z)$ . Thus, we obtain
 the
following statement:

\begin{lemma}
{\sl If the function $w_{\sigma ,c}$, defined by {\rm(\ref{3.1})},
  satisfies  the interpolation
conditions {\rm(\ref{3.12})}, then the {\rm FMI} {\rm (\ref{1.5})}
 (with $w$ replaced by $w_{\sigma ,c}$)
 is satisfied for every $z\in \Bbb C \setminus \Bbb T$ , where
$A$ is defined by {\rm(\ref{2.4})}, $B_w$ is defined by {\rm(\ref{1.8})},
{\rm (\ref{2.2})}, {\rm (\ref{2.3})} and
 $C_w$ is defined by \rm{(\ref{1.6})}}.
\end{lemma}

According to the Riesz-Herglotz theorem,
 each function $w$ from the considered class
has the representation $w=w_{\sigma,c}$. Thus, the following result holds:

\begin{theo} ({\cmss From interpolation conditions to FMI($\cal NP$)}).
{\sl Let interpolation data for $\cal NP$ problem be given.
 Let $w$ be a function, which belongs to the class \Cgot\,($\Bbb D$).
 If the function $w$ satisfies the
 interpolation conditions {\rm{(\ref{2.1})}}, then the {\rm FMI($\cal NP$)} for
 this function
 (with $A$ and $ v$ constructed from the given interpolation data)
 is satisfied for every $z\in \Bbb C \setminus \Bbb T$.
}
\end{theo}

We have stated this (well known) derivation of the FMI ($\cal NP$) because the
formulas (\ref{3.4}) and (\ref{3.6}) are a very convenient starting point
to guess formulas for transformations of FMI.

\vspace*{0.25cm} 
\begin{minipage}{15.0cm}
\section{\hspace{-0.4cm}.\hspace{0.2cm} THE HAMBURGER MOMENT PROBLEM \newline
AS A CLASSICAL INTERPOLATION PROBLEM}
\end{minipage}\\[-0.2cm]
\setcounter{equation}{0}

This problem can be considered as a classical interpolation problem
in the class \Rgot\,($\Bbb H$).

FORMULATION OF THE HAMBURGER MOMENT PROBLEM.
{\sl  {\cmss The data} of the Hamburger problem is a finite sequence
$s_0,s_1, \ldots, s_{2n-1},s_{2n}$ of real numbers.
A nonnegative measure $\sigma$ on the real numbers is said to be
{\cmss a solution of the Hamburger moment problem} (with these data),
if its power moments
\begin{equation}
s_k(\sigma)=\int_{\Bbb R} \lambda^k d\sigma(\lambda)
\qquad (\,  k=0,1,\ldots,2n-1,2n \, ) 
\label{4.4}
\end{equation}
exist and satisfy {\cmss the moment conditions}
\begin{equation}
\mbox{\rm i).\ \ }
s_k(\sigma)=s_k
\qquad (\,  k=0,1,\ldots,2n-1\, )\, ;
\qquad
\mbox{\rm ii).\ \ }
 s_{2n}(\sigma) \leq s_{2n}. 
\label{4.5}
\end{equation}
Measures $\sigma$ satisfying these moment conditions are sought.}
\vspace*{0.3cm}

At first glance the formulated moment problem does not look like an
interpolation problem. However, this problem can be reformulated
as a classical interpolation problem.

Namely, let $\sigma$ be a nonnegative measure on $\Bbb R$ which is finite:
$s_0 (\sigma)< \infty  .$
We associate with this measure $\sigma$ the function $w_{\sigma}:$
\begin{equation}
w_{\sigma}(z)= \int_{\Bbb R}\frac{d\sigma (\lambda)}{\lambda - z}
\qquad (\,  
z \in {\Bbb C}\setminus{\Bbb R}\, ) 
\label{4.6}
\end{equation}
This function $w_{\sigma}$ belongs to the class $\mbox{\Rgot}_0$($\Bbb H$).

The following result is a version of the Riesz - Herglotz theorem 
for the upper half-plane.

THEOREM ({\cmss Nevanlinna}).
{\sl Let $w$ be a function from the class
$\mbox{\Rgot}_0$\,($\Bbb H$). Then this function $w$ is representable
in the form  (\ref{4.6}), with some finite nonnegative 
measure \mbox{$\sigma:\, \sigma \geq 0,$} $s_0(\sigma)< \infty.$
This measure $\sigma$ is determined from the function $w$ uniquely.}

It turns out that if a measure $\sigma$ solves the Hamburger moment problem
 (\ref{4.5}), then the function $w_{\sigma},$ associated with this
measure $\sigma,$ satisfies some asymptotic relation. To obtain such a
 relation, we consider the functions $w_{\sigma, k}:$
\begin{equation}
w_{\sigma, k}(z)=
 \int_{\Bbb R}{\lambda^k \frac{d\sigma (\lambda)}{\lambda - z}}
\qquad (\,  k=0,1,2,\ldots ,2n
\, ). 
\label{4.7}
\end{equation}
(In this notation, $w_{\sigma}=w_{\sigma,0}$).
 Assume that a measure $\sigma \geq 0$
on $\Bbb R$ has the moment $s_{2n}(\sigma)$ (and hence, also the moments
$s_0(\sigma), \ldots, s_{2n-1}(\sigma)$). Integrating the identity
\begin{equation}
\frac{\lambda^k}{\lambda -z}=\frac{z^k}{\lambda -z} +
\sum_{0 \leq j \leq k-1}z^{k-1-j}\lambda^j
\label{4.8}
\end{equation}
with respect to the measure $\sigma$,  we come to the equality
\begin{equation}
w_{\sigma,k}(z)=z^k \left ( w_{\sigma}(z) + \sum_{0 \leq j \leq k-1}
\frac{s_j(\sigma)}{z^{j+1}} \right )
\qquad (\,  k=0,1,2,\ldots ,2n
\, ). 
\label{4.9}
\end{equation}
Since
\begin{equation}
w_{\sigma,2n}(z)=-\frac{s_{2n}(\sigma)}{z}(1+o(1))
\qquad (\, |z| \to \infty,\ z=iy
\, ), 
\label{4.10}
\end{equation}
it follows from (\ref{4.9}) (with $k=2n$) that
\begin{equation}
z^{2n} \left ( w_{\sigma}(z) + \sum_{0 \leq j \leq 2n-1}
\frac{s_j(\sigma)}{z^{j+1}} \right )=-\frac{s_{2n}(\sigma)}{z}(1+o(1))
 \qquad (\, |z| \to \infty,\ z=iy
\, ). 
\label{4.11}
\end{equation}

The asymptotic relation (\ref{4.11}), together with (\ref{4.5}),(\ref{4.9})
suggests  the following:

{\sl Given the function $w$ of the class \Rgot\,($H$)
 and a set of real
numbers $s_0,s_1, \ldots , s_{2n-1},$ it has to be profitable to consider
the functions $b_{w,k}(z)=b_{w,k}(z;s_0,s_1, \ldots ,s_{k-1}) :$
\begin{equation}
b_{w,k}(z)=z^kw(z)+\sum_{0 \leq j \leq k-1}z^{k-1-j}s_j
\qquad (\,  k=0,1,2,\ldots ,2n
\, ) 
\label{4.12}
\end{equation}
and the asymptotic relation of the form}
\begin{equation}
|b_{w,k}(z)|=O(|z|^{-1})
\qquad (\, |z| \to \infty,\ z=iy
\, ). 
\label{4.13}
\end{equation}
In this notation the equality (\ref{4.9}) means that 
\begin{equation}
w_{\sigma,k}(z)=b_{w_{\sigma},k}(z;s_0(\sigma), \ldots , s_{k-1}(\sigma))
\label{4.14}
\end{equation}
From (\ref{4.11}) and (\ref{4.14}) it follows that:

{\sl If a measure $\sigma \geq 0$ on $\Bbb R$ satisfies the moment conditions
(\ref{4.5}), then the asymptotic relation
\begin{equation}
|b_{w_{\sigma},2n}(z;s_0, \ldots , s_{2n-1})| \leq \frac{s_{2n}}{|z|}
(1+o(1))
\qquad (\, |z| \to \infty,\ z=iy
\, ) 
\label{4.15}
\end{equation}
holds.}

It is remarkable that the last statement can be inverted.

THEOREM ({\cmss Hamburger}).
{\sl Let $w$ be a function which belongs to the class \Rgot\,($\Bbb H$)
 and let $s_0,s_1, \ldots , s_{2n-1}$ be real numbers.
Assume that the function $w$ satisfies the asymptotic condition
\begin{equation}
|b_{w,2n}(z;s_0, \ldots , s_{2n-1})| = O(|z|^{-1})
\qquad (\, |z| \to \infty,\ z=iy
\, ) 
\label{4.16}
\end{equation}
(where $b_{w,2n}$ is defined in (\ref{4.12})). Then the function $w$
has the representation of the form (\ref{4.6}), with
a nonnegative measure $\sigma,$ which has $2n$-th moment: $s_{2n}(\sigma)
< \infty.$ Moreover, 
\begin{equation}
s_0(\sigma)=s_0,s_1(\sigma)=s_1, \ldots, s_{2n-1}(\sigma)=
s_{2n-1},
\label{4.17}
\end{equation}
\begin{equation}
s_{2n}(\sigma)=\lim_{{|z| \to \infty} \atop {z=iy}} (-z)b_{w,2n}
(z;s_0,s_1, \ldots, s_{2n-1})
\label{4.18}
\end{equation}
}

This theorem was proved by Hamburger (\cite{H}, Theorem $IX$).
It is reproduced in the monograph by N. Akhiezer (\cite{A}, Theorem 2.3.1).
The proof which was presented by Hamburger is based on
a ``step by step'' algorithm. Another proof of this theorem, and
its far reaching generalizations, is presented in \cite{K1}.

Thus the Hamburger moment problem can be reformulated as the following
interpolation problem: \\[0.3cm]
\hspace*{0.5cm}
\begin{minipage}{15.2cm}
{\cmss Function class}: {\sl the class \Rgot\,($\Bbb H$). \\
{\cmss Interpolation data}: a finite sequence $s_0,s_1, \ldots,s_{2n}$ of
real numbers.  \\
The asymptotic relation}
\end{minipage}
\begin{equation}
\left |z^{2n} \left ( w(z) + \sum_{0 \leq j \leq 2n-1}
\frac{s_j}{z^{j+1}} \right ) \right | \leq \frac{s_{2n}}{|z|}(1+o(1))
 \qquad (\, |z| \to \infty,\ z=iy
\, ) 
\label{4.19}
\end{equation}
\hspace*{0.5cm}
\begin{minipage}{15.2cm}
{\sl is considered as an {\cmss interpolation condition}.
(The point $z=\infty$ is a multiple interpolation node which
lies on the boundary of the upper half-plane ${\Bbb H}.$ Its
multiplicity  equals  $2n$). We seek functions $w$ from this 
class which satisfy the  condition} (\ref{4.19}).
\end{minipage}\footnote
{\  Strictly speaking, the considered problem has two
 interpolation nodes
which are symmetric with respect to the real axis and are located
at the points $+i\cdot \infty$ and $-i\cdot \infty.$ The multiplicity of each
 of them equals  $n.$}

\begin{rem}
 $i$). Assume that a function $w$ from the class \Rgot\,($\Bbb H$) 
satisfies the condition (\ref{4.16}). Suppose that we also know
 (for example, from the Hamburger theorem), that 
$w=w_{\sigma}$, where $s_{2n}(\sigma)<\infty$.
Then we can construct the function $w_{\sigma , 2n}$  by (\ref{4.9}).
Comparing the asymptotics (\ref{4.16}) and (\ref{4.10}), we conclude, that 
$b_{w, 2n}=w_{\sigma , 2n}$. Hence, the moment condition
(\ref{4.5}.\hspace{1pt}i) is satisfied, as well as the condition
\begin{equation}
\left |z^{2n} \left ( w(z) + \sum_{0 \leq j \leq 2n-1}
\frac{s_j}{z^{j+1}} \right ) \right | \leq \frac{s_{2n}(\sigma)}{|z|}(1+o(1))
 \qquad (\, |z| \to \infty,\ z=iy
\, ). 
\label{4.20}
\end{equation}
Moreover, the function $b_{w,2n}(z;s_0,\ldots, s_{2n-1})$ belongs
to the class $\mbox{\Rgot}_0$($\Bbb H$).
(If $d\sigma (\lambda)$ is a measure which represents $w,$, then the measure
$\lambda^{2n}d\sigma (\lambda)$ represents the function $b_{w,2n}$).  

$ii$). Assume now  that  the function
 $b_{w,2n}(z;s_0,\ldots, s_{2n-1})$ belongs to the class 
$\mbox{\Rgot}_0$($\Bbb H$).
Then, by the Nevanlinna'sn  theorem, the function $b_{w,2n}$ has the form
 $w_{\tau}$ for some
$d\tau \geq 0,s_0(\tau)<\infty.$
Thus,
$$
\int_{\Bbb R}\frac{d\tau (\lambda)}{\lambda - z}=
z^{2n}\int_{\Bbb R}\frac{d\sigma (\lambda)}{\lambda - z}+
\sum_{0 \leq j \leq 2n-1}s_jz^{2n-1-j}
$$
Applying the generalized Stieltjes inversion formula (\cite{KaKr},$\S 2$),
 we conclude that
$d\tau (\lambda)={\lambda}^{2n}d\sigma (\lambda).$
Hence, $\int_{\Bbb R}{\lambda}^{2n}d\sigma (\lambda)=
\int_{\Bbb R}d\tau (\lambda)< \infty.$ Thus, 
$b_{w,2n}=w_{\sigma,2n }$, and (\ref{4.20}) is satisfied.
\end{rem}

\vspace*{0.2cm}
\begin{minipage}{15.0cm}
\section{\hspace{-0.4cm}.\hspace{0.2cm}DERIVATION  OF THE
 FMI\,\boldmath$(\, \cal H\,)$\unboldmath }
\end{minipage}\\[-0.2cm]
\setcounter{equation}{0}

Given the Hamburger moment problem with data $s_0, s_1,\, \ldots,\, s_{2n}$,
we associate with this problem {\sl the Pick matrix}
\begin{equation} 
A=
\left[
\begin{array}{cccc}
s_0&s_1&\cdots &s_n\cr
s_1&s_2&\cdots &s_{n+1}\cr
\cdots&\cdots&\cdots&\cdots\cr
s_{n-1}&s_{n}&\cdots &s_{2n-1}\cr
s_{n}&s_{n+1}&\cdots &s_{2n}\cr
\end{array}
\right],
\label{5.1}
\end{equation}
and the vectors of the interpolation data
\begin{equation} 
u=
\left[
\begin{array}{c}
1\cr
0\cr
\vdots\cr
0\cr
0\cr
\end{array}
\right]
\qquad
\mbox{and}
\qquad
v=
\left[
\begin{array}{c}
0\cr
-s_{0}\cr
\vdots\cr
-s_{n-2}\cr
-s_{n-1}\cr
\end{array}
\right]\cdot
\label{5.2}
\end{equation}
The matrix, which is responsible for interpolation knots (with multiplicity) is:
\begin{equation} 
T=
\left.
\left[
\begin{array}{cccccc}
0&0&\cdots&0&0&0\cr
1&0&\cdots&0&0&0\cr
0&1&\cdots&0&0&0\cr
\cdots&\cdots&\cdots&\cdots&\cdots&\cdots\cr
0&0&\cdots&0&0&0\cr
0&0&\cdots&1&0&0\cr
0&0&\cdots&0&1&0\cr
\end{array}
\right]
\right\}{\scriptstyle (n+1)}\ \cdot
\label{5.3}
\end{equation}
The Fundamental Identity (\ref{1.11}) for this chois of
 $u,\, v,\, T$ and $A$ can be checked straightforwardly.

Now we derive th Fundamental Matrix Inequality for the Hamburger Moment Problem
(\,FMI\,($\cal H$)\,). Let $d\sigma (\lambda)$ be a nonnegative measure on 
$\Bbb R$ for which the $2n$ th moment is finite: $s_{2n}(\sigma)< \infty$. 
The following inequality is clear:
\begin{equation}
\raisebox{0.65cm}{$\Int$ \hspace{-0.55cm}
\raisebox{-1.4cm}{$\scriptstyle\Bbb R$}}\hspace{0.25cm}
\left[
\begin{array}{c}
\vspace*{-0.1cm}
(I-\lambda T)^{-1}u \cr
\vspace*{-0.15cm}
\cdot - \cdot - \cdot - \cdot 
\vspace*{0.15cm}\cr 
\displaystyle (\bar{\lambda}-\bar{z})^{-1}\cr 
\end{array}
\right]\cdot
d\sigma(\lambda)\cdot
\left[
\begin{array}{ccc}
u^{\ast}(I-\bar{\lambda}T^{\star})^{-1}&\hspace{-0.2cm} 
\begin{array}{c}
|\vspace*{-0.1cm}\cr
\vspace*{-0.1cm}
\cdot
 \cr
\vspace*{-0.1cm}|\vspace*{0.1cm}\cr
\end{array}
\hspace{-0.2cm}  &\displaystyle (\lambda - z)^{-1}
\end{array}
\right]\,
\geq 0 .
\label{5.4}
\end{equation}
This inequality has the form
\begin{equation}
\left[
\begin{array}{cc}
A_{\sigma}&B_{\sigma}(z)\cr
 \vspace*{-0.1cm}&\vspace*{-0.1cm}\cr
B^{\ast}_{\sigma}(z)&
\displaystyle\frac{w_{\sigma}(z)-w_{\sigma}^{\ast}(z)}{z-\bar{z}}\cr
\end{array}
\right]\geq 0,
\label{5.5}
\end{equation}
where the function $w_{\sigma}$ is defined by (\ref{4.6}).
It is clear that
\begin{equation}
A_{\sigma}=\int\limits_{\Bbb R}
(I-\lambda T)^{-1}u \cdot d\sigma (\lambda)
\cdot
u^{\ast}(I-\lambda T^{\ast})^{-1},
\label{5.6}
\end{equation}
where
\begin{equation}
A_{\sigma}=\|a_{\sigma ,kl}\|_{0\leq k.l \leq n},
\qquad a_{\sigma ,kl}=s_{k+l}(\sigma)\quad (0\leq k,\, l \leq n).
\label{5.7}
\end{equation}
It is also clear, that
\begin{equation}
B_{\sigma}(z)=\int\limits_{\Bbb R}\frac{(I-\lambda T)^{-1}u}{\lambda -z}
\: d\sigma(\lambda ). 
\label{5.8}
\end{equation}
Since
\begin{equation}
\frac{(I-\lambda T)^{-1}}{\lambda-z}= (I-zT)^{-1}
\left(
\frac{1}{\lambda - z} + T(I-\lambda T)^{-1}
\right),
\label{5.9}
\end{equation}
it follows that
\begin{equation}
B_{\sigma}(z)=(I-zT)^{-1}\big(u\cdot w_{\sigma}(z)-v_{\sigma}\big),
\label{5.10}
\end{equation}
where
\begin{equation}
v_{\sigma}=-\int\limits_{\Bbb R}T(I-\lambda T)^{-1}u\,d\sigma (\lambda ).
\label{5.11}
\end{equation}
From the concrete expressions (\ref{5.2}) and  (\ref{5.3}) 
 for $u$ and $T$ it is not difficult to see that
\begin{equation} 
v_{\sigma}=
\left[
\begin{array}{c}
0\cr
-s_{0}(\sigma)\cr
\vdots\cr
-s_{n-2}(\sigma)\cr
-s_{n-1}(\sigma)\cr
\end{array}
\right].
\label{5.12}
\end{equation}

Assume now,that the measure $\sigma$ satisfies the moment conditions 
(\ref{4.5}). Then, according to (\ref{5.2}) and (\ref{5.12}),
 $v_{\sigma}=v$, and according to (\ref{5.1}) and (\ref{5.7}),
$a_{\sigma,kl}=a_{kl}\ \ (0\leq k+l<2n,\, a_{\sigma,nn}\leq a_{nn}$, hence,
$A_{\sigma}\leq A$. Thus, we obtain

\begin{theo} {\cmss(From the moment conditions to the FMI\,($\cal H$))}.
{\sl  Let interpolation data for the Hamburger moment problem be given.
Let $w$ be a function of the form  {\rm (\ref{4.6})}, where the 
measure $\sigma$ satisfies the moment conditions {\rm (\ref{4.5})} (or,
what is the same according to Hamburger, the interpolation condition
{\rm (\ref{4.19})} is satisfied). Then the
 {\rm FMI($\cal H$)} {\rm (\ref{1.5})} holds for
this function $w$ at every point $z\in {\Bbb C} \setminus {\Bbb R}$,
where $A$ is defined by {\rm (\ref{5.1})},
 $C_{w}$ is defined by {\rm (\ref{1.7})} and
$B_{w}$ is defined by (\ref{1.9}), (\ref{5.2}), (\ref{5.3}).
}
\end{theo}
\newpage
\vspace*{0.2cm}
\begin{minipage}{15.0cm}
\section{\hspace{-0.4cm}.\hspace{0.2cm}TRANSFORMATION OF THE FMI
 \boldmath$(\, \cal H\,)$\unboldmath }
\end{minipage}\\[-0.2cm]
\setcounter{equation}{0}

Let $s_{0},\, \ldots,\, s_{2n}$ be interpolation data for the Hamburger
 moment problem. Then the Pick matrix A is defined by (\ref{5.1}),
the interpolation nodes matrix $T$ be defined by (\ref{5.3}) and interpolation
values vectors $u$ and $v$ are defined by (\ref{5.2}). Given a function $w$,
which is holomorphic
 in ${\Bbb C} \setminus {\Bbb R}$ and satisfies the symmetry
conditions (\ref{1.3}), assume that the FMI\,($\cal H$\,)
\begin{equation}
\left[
\begin{array}{ccc}
A&\hspace{-0.1cm} | \hspace{-0.1cm} &B_{w}(z)\cr
-\cdot  - &\hspace{-0.1cm}\cdot \hspace{-0.1cm} & -\cdot  -\cdot  -\cdot  - \cr
B^{\ast}_{w}(z)&\hspace{-0.1cm} | \hspace{-0.1cm}
 &\displaystyle\frac{w(z)-w^{\ast}(z)}{z-\bar{z}}\cr
\end{array}
\right]
\geq 0
\label{6.1}
\end{equation}
is satisfied for every $z\in {\Bbb C} \setminus {\Bbb R}$. Here $B_{w}$ is
defined by (\ref{1.9}), (\ref{5.2}), (\ref{5.3}),
 or in detail,
\begin{equation}
B_{w}(z)=
\left[
\begin{array}{c}
0\cr
b_{w,0}(z)\cr
\cdots \cr
b_{w,n-1}(z)\cr
\end{array}
\right] \cdot
\label{6.2}
\end{equation}
Our goal is to extract  interpolation information from this FMI.
Of course, from (\ref{6.1}) it follows, that the function $w$ satisfies
the positivity condition (\ref{1.4}). Proceeding in the same way, as in
the Proof of Theorem 2.1, we have to consider the
``subinequalities'' (\ref{2.7})
of the inequality (\ref{6.1}). The most  information  which we can
obtain in this way from (\ref{6.1})  is contained in the subinequality
\begin{equation}
\left[
\begin{array}{cc}
s_{2n}&b_{w,n-1}\cr
\vspace{-0.1cm} &\vspace{-0.1cm} \cr
b^{\ast}_{w,n-1}&\frac{w(z)-w^{\ast}(z)}{z-\bar{z}}\cr
\end{array}
\right]
\geq 0.
\label{6.3}
\end{equation}
First and foremost, from  (\ref{6.3}) we obtain the estimate (\ref{1.4.a})
 for $w$. By the Nevanlinna Theorem, the function $w$ has the form $w_{\sigma}$
for some nonnegative measure $\sigma$ with  $s_{0}(\sigma)< \infty$. Moreover,
 the estimate 
$
|b_{w,n-1}(iy)|=O(|y|^{-1})
$
as
$
y  \uparrow \infty
$
follows from (\ref{6.3}).
This is not enough since the function $b_{w,n-1}$ contains the interpolation
data $s_{0},\,s_{1},\, \ldots,\, s_{n-1}$ only, and does not contain
the data  $s_{n},\,s_{n+1},\, \ldots,\, s_{2n-1}$ at all.
 We need to obtain the condition (\ref{4.19}) from (\ref{6.1}).
Clearly, it is impossible to extract the condition (\ref{4.19}) by considering
``subinequalities'' of the inequality (\ref{6.1}). More generally, it is
 impossible to obtain (\ref{4.19}) from any inequality of the form
 (\ref{1.12}) when the framing matrix $M$ does not depend on $z$ because
the data $s_{n},\,s_{n+1},\, \ldots,\, s_{2n-1},\,s_{2n} $ appear in the 
block $A$ only, which does not depend on $z$.

Therefore, in order  to extract (\ref{4.19}) from (\ref{6.1})
(if it is at all possible),
we have to choose a matrix $M$ in (\ref{1.12}), which depends on $z$.
 To understand how to do this we return to the derivation of the
 FMI\,($\cal H$)\,. Let us consider the inequality (\ref{5.5}). It contains the
functions $w_{\sigma ,k}=b_{w_{\sigma ,k}}$
 with $k=0, 1, \ldots , n-1$ only. However, we need
the function $w_{\sigma, 2n-1}$. The only information which is available for us
is the block $A_{\sigma}$, which is defined by
 (\ref{5.6}) and (\ref{5.7}). The Hankel matrix
$A_{\sigma}$ is related to the Hankel matrix
\begin{equation}
W_{\sigma}(z)=\|w_{\sigma, kl}\,(z)\|_{0\leq k,l \leq n}\, ,
\label{6.4}
\end{equation}
with entries
\begin{equation}
w_{\sigma, kl}\,(z)=\int\limits_{\Bbb R}\,
\lambda^k\cdot\frac{d\sigma (\lambda)}{\lambda - z}\cdot\lambda^l 
\qquad (0\leq k,l \leq n).
\label{6.5}
\end{equation}
$k,l$-entries of the matrix $W_{\sigma}$ with $k+l < n$ are the same
functions which appear  in the column $B_{\sigma}$. The entries with
$n\leq k+l \leq 2n$ are exactly those which we need.
Thus, the problem is  to obtain the matrix $W_{\sigma}$ from the 
matrix $A_{\sigma}$. According to (\ref{6.5}), (\ref{5.2}) and (\ref{5.3}),
\begin{equation}
W_{\sigma}(z)=\int\limits_{\Bbb R}
(I-\lambda T)^{-1}u\cdot
\frac{d\sigma (\lambda)}{\lambda - z}
\cdot u^{\ast}
(I-\lambda T^{\ast}). 
\label{6.6}
 \end{equation}
Comparing  (\ref{6.6}) with  (\ref{5.6}) we see that we
 have to replace
$(I-\lambda T)^{-1}$ with
 $\displaystyle\frac{(I-\lambda T)^{-1}}{\lambda - z}$\  in (\ref{5.6}).
Let us turn to the identity (\ref{5.9}):
\begin{equation}
T(I-zT)^{-1}\cdot(I-\lambda T)^{-1}u=
\frac{(I-\lambda T)^{-1}}{\lambda - z}\, u - 
\frac{(I-zT)^{-1}}{\lambda - z}\, u\,  .
\label{6.7}
\end{equation}
From (\ref{6.6}) and (\ref{6.7}) it follows that
\begin{equation}
T(I-zT)^{-1}A_{\sigma}=W_{\sigma}(z)-(I-zT)^{-1}u\cdot
\int\limits_{\Bbb R}\, d\sigma(\lambda)\,
\frac{u^{\ast}(I-\lambda T^{\ast})^{-1}}{\lambda - z}\cdot
\label{6.8}
\end{equation} 
Taking into account (\ref{5.8}), we obtain the equality
\begin{equation}
W_{\sigma}(z)=T((I-zT)^{-1})\,A_{\sigma}+(I-zT)^{-1}u
\cdot B^{\ast}_{\sigma}(\bar{z}). 
\label{6.9}
\end{equation}
The equality (\ref{6.9}) provide us a heuristic reason for the following

\begin{defi}
{\sl 
Given a Hermitian matrix $A$, a matrix $T$ and vectors $u$ and $v$, which
satisfy the Fundamental Identity {\rm (\ref{1.11})}, we associate with each
 function  $w$, which is holomorphic in
 $\Bbb C \setminus \Bbb R$ and satisfies the symmetry condition
{\rm (\ref{1.3})},
the function $W_{w}$:
\begin{equation}
W_{w}(z)=T(I-zT)^{-1}\,A +(I-zT)^{-1}u
\cdot B^{\ast}_{w}(\bar{z}). 
\label{6.10}
\end{equation} 
or, in detail,
\begin{eqnarray}
W_{w}(z)=T(I-zT)^{-1} \,A - 
(I-zT)^{-1}u \cdot v^{\ast}(I-zT^{\ast})^{-1}  \cr
     \vspace{-0.3cm}                                     \cr
+\, (I-zT)^{-1}u\cdot w(z) \cdot  u^{\ast}(I-zT^{\ast})^{-1}.
\label{6.11}
\end{eqnarray}
}
\end{defi}

\begin{lemma}
{\sl
The matrix function $W_{w}$ satisfies the same symmetry condition as
 that the function $w$:
\vspace{-0.2cm}
\begin{equation}
W_{w}(z)=W^{\ast}_{w}(\bar{z}) \qquad (z\in \Bbb C \setminus \Bbb R).
\label{6.12}
\end{equation}
}
\end{lemma}

Straightforward calculation gives us the explicit expression for $W_{w}(z)$:
\begin{equation}
W_{w}(z)=\|b_{w ,k+l}\,(z)\|_{0\leq k,l \leq n}
\label{6.13}
\end{equation}
Thus, the matrix-function $W_{w}$ is exactly what we need:  it contains 
the function $b_{w, 2n}$.
 In particular,
 from the formula it follows that the matrix $W_{w}(z)$ is a Hankel
matrix. However, the Hankel structure of the matrix $W_{w}(z)$ can be
obtained in a less special way, i.e. by using the FI
 (\ref{1.11}) only:
\vspace*{0.2cm}

\begin{lemma}
{\sl
The matrix $W_{w}(z)$ satisfies the following identity\footnote
{\  The equality (\ref{6.14}),
 considered as an equation with respect to the matrix $W_{w}(z)$,  can be
used to calculate this matrix.}:
\vspace*{-0.2cm}
\begin{equation}
T\,W_{w}(z)-W_{w}(z)\,T^{\ast}
=u\cdot\varphi^{\ast}_{w}(z)-\varphi_{w}(\bar{z})\cdot u^{\ast},
\quad \mbox{\sl where}\quad
\varphi_{w}(z)= -T(I-zT)^{-1}\Big( u\cdot w(z) - v \Big).
\label{6.14}
\end{equation}
}
\end{lemma}

\begin{lemma}
{\sl
For the Hamburger moment problem,
 the function $w(z)$ and the column $B_{w}(z)$ can be recovered from the
matrix-function $W_{w}(z)$ in the following way:
\begin{equation}
w(z)=e_0\cdot W_{w}(z)\cdot e^{\ast}_0\,, \qquad 
 B_{w}(z)=W_{w}(z)\cdot e^{\ast}_0\, ,
\label{6.15}
\end{equation}
where
$
e_0=
\left[
\begin{array}{cccc}
1&0&\cdots &0\cr
\end{array}
\right]
$
 is a $(n+1)\times 1$ vector.
}
\end{lemma}

PROOF. The formulas in (\ref{6.15}) follows from the equalities
\begin{equation}
e_0\, T=0\, ,\quad  e_0\, u= 1 \quad  \mbox{\sl and }\quad  e_0\, v= 0.
\label{6.17}
\end{equation}
 \hfill \raisebox{1.0cm}{$\boxx$}
\vspace*{-0.6cm}

\begin{rem}
{\sl
The proof of the lemma depends on the
equalities (\ref{6.17}), not on the FI (\ref{1.11}).  It is specific for the  problem in question.
}
\end{rem}
\vspace*{0.2cm}

Let us turn to the FMI (\ref{6.1}). It is clear that the matrix
 $W_{w}(\bar{z})$ appears in the product

\begin{equation}
\left[
\begin{array}{ccc}
T(I-\bar{z}T)^{-1}&\vdots&(I-\bar{z}T)^{-1}u\cr
\end{array}
\right]
\cdot
\left[
\begin{array}{ccc}
A&\hspace{-0.1cm} | \hspace{-0.1cm} &B_{w}(z)\cr
-\cdot  - &\hspace{-0.1cm}\cdot \hspace{-0.1cm} & -\cdot  -\cdot  -\cdot  - \cr
B^{\ast}_{w}(z)&\hspace{-0.1cm} | \hspace{-0.1cm}
 &\displaystyle\frac{w(z)-w^{\ast}(z)}{z-\bar{z}}\cr
\end{array}
\right]
\label{6.18}
\end{equation}
In order to to transform the FMI (\ref{6.1}), we have to ``frame''
 it according to
 (\ref{1.12}), where now the matrix $M$ depends on $z$. It is clear that
the row
$
\left[
\begin{array}{ccc}
T(I-\bar{z}T^{-1})&\vdots&(I-\bar{z}T)^{-1}u\cr
\end{array}
\right]
$
ought to be one of the rows of the matrix $M(z)$.  There are two main
possibilities. Either the mentioned row is the first row of the matrix $M$:
\begin{equation}
M_{1}(z)=
\left[
\begin{array}{ccc}
T(I-\bar{z}T)^{-1}&\ &(I-\bar{z}T)^{-1}u \cr
\vspace{-0.1cm}& \vspace{-0.1cm}& \vspace{-0.1cm}\cr      
0                   &\ &   1 \cr
\end{array}
\right]\, ,
\label{6.19}
\end{equation}
or the mentioned row is the second row of the matrix $M$:
\begin{equation}
M_{2}(z)=
\left[
\begin{array}{ccc}
I                   &\ &   0 \cr
\vspace{-0.1cm}& \vspace{-0.1cm}& \vspace{-0.1cm}\cr      
T(I-\bar{z}T)^{-1}&\ &(I-\bar{z}T)^{-1}u\cr
\end{array}
\right]\cdot
\label{6.20}
\end{equation}
Upon performing the matrix multiplications, we obtain
 (after some calculations with
the matrix entries):
\begin{equation}
M_{1}(z)\cdot
\left[
\begin{array}{ccc}
A &B_{w}(z)\cr
& \cr
B^{\ast}_{w}(z)&\displaystyle\frac{w(z)-w^{\ast}(z)}{z-\bar{z}}\cr
\end{array}
\right]
\cdot
M^{\ast}_{1}(z)=
\left[
\begin{array}{ccc}
\displaystyle\frac{W_{w}(z)-W^{\ast}_{w}(z)}{z-\bar{z}}& 
  &\displaystyle\frac{B_{w}(z)-B_{w}(\bar{z})}{z-\bar{z}}\cr
&  & \cr
\displaystyle\frac{B^{\ast}_{w}(\bar{z})-B^{\ast}_{w}(z)}{z-\bar{z}}&
 &\displaystyle\frac{w(z)-w^{\ast}(z)}{z-\bar{z}} \cr
\end{array}
\right]
\label{6.21}
\end{equation}
and
\begin{equation}
M_{2}(z)\cdot
\left[
\begin{array}{ccc}
A &B_{w}(z)\cr
& \cr
B^{\ast}_{w}(z)&\displaystyle\frac{w(z)-w^{\ast}(z)}{z-\bar{z}}\cr
\end{array}
\right]
\cdot
M^{\ast}_{2}(z)=
\left[
\begin{array}{ccc}
\displaystyle A& 
  &\displaystyle W_{w}(z)\cr
       & & \cr

W^{\ast}_{w}(z)  & &
\displaystyle \frac{W_{w}(z)-W^{\ast}_{w}(z)}{z-\bar{z}}\cr
\end{array}
\right]\cdot
\label{6.22}
\end{equation}
The  calculations with the matrix entries are based essentially
on the following consequence of the FI (\ref{1.11}):
\vspace*{0.2cm}

\begin{lemma}
{\sl
The identity
\vspace*{-0.2cm}
\begin{eqnarray}
T(I-zT)^{-1}\cdot A\cdot( I-\bar{z}T^{\ast})^{-1}T^{\ast}=
\hspace{7.5cm}\cr
\vspace{-0.2cm}\cr
=\frac{T(I-zT)^{-1}\, A - A\, (I-\bar{z}T^{\ast})^{-1}T^{\ast}}{z-\bar{z}}
-(I-zT)^{-1}\cdot
 \frac{uv^{\ast}-vu^{\ast}}{z-\bar{z}}
\cdot(I-\bar{z}T^{\ast})^{-1} 
\label{6.23}
\end{eqnarray}
\vspace*{-0.2cm}
holds.
}
\end{lemma}
\vspace*{0.8cm}

\begin{minipage}{15.0cm}
\section{\hspace{-0.4cm}.\hspace{0.2cm}USING OF THE
 TFMI\,\boldmath$( \cal H )$\unboldmath\  -- \
 \newline
 FROM THE FMI\,\boldmath$(\, \cal H\,)$\unboldmath\ TO INTERPOLATION
INFORMATION}
\end{minipage}\\[-0.2cm]
\setcounter{equation}{0}

We consider two kinds of Transformed Fundamental Matrix Inequalities (for the
Hamburger problem):
 $\mbox{\rm TFMI}_{\mbox{\scriptsize\rm I}}(\cal H)$ and
 $\mbox{\rm TFMI }_{\mbox{\scriptsize\rm II}}(\cal H)$.
 
The $\mbox{\rm TFMI}_{\mbox{\scriptsize\rm I}}(\cal H)$ is of the form
\begin{equation}
\left[
\begin{array}{ccc}
\displaystyle\frac{W_{w}(z)-W^{\ast}_{w}(z)}{z-\bar{z}}& 
  &\displaystyle\frac{B_{w}(z)-B_{w}(\bar{z})}{z-\bar{z}}\cr
&  & \cr
\displaystyle\frac{B^{\ast}_{w}(\bar{z})-B^{\ast}_{w}(z)}{z-\bar{z}}&
 &\displaystyle\frac{w(z)-w^{\ast}(z)}{z-\bar{z}} \cr
\end{array}
\right]
\geq 0.
\label{7.1}
\end{equation}
The $\mbox{\rm TFMI}_{\mbox{\scriptsize\rm II}}(\cal H)$ is of the form
\begin{equation}
\left[
\begin{array}{ccc}
\displaystyle A& 
  &\displaystyle W_{w}(z)\cr
       & & \cr

W^{\ast}_{w}(z)  & &
\displaystyle \frac{W_{w}(z)-W^{\ast}_{w}(z)}{z-\bar{z}}\cr
\end{array}
\right]
\geq 0.
\label{7.2}
\end{equation}
We see that both of the TFMI's contain the function $W_{w}(z)$. Now the problem
of  extracting  interpolation information from the TFMI arises.

Now we will discuss the extent to which
 the FMI\,($\cal H$)\, and the TFMI\,(\,$\cal H$)\,
 are equivalent.
In view of (\ref{6.21}) and (\ref{6.22}), it is clear that
\begin{equation}
\mbox{FMI}\,(\cal H)\, \ \ \Rightarrow \ \ 
\mbox{\rm TFMI}_{\mbox{\scriptsize\rm I}}\,(\cal H)
\label{7.3}
\end{equation}
and
\begin{equation}
\mbox{FMI}\,(\cal H)\, \ \ \Rightarrow \ \ 
\mbox{\rm TFMI}_{\mbox{\scriptsize\rm II}}\,(\cal H).
\label{7.4}
\end{equation}
More formally:

\begin{lemma}
 {\sl If the {\rm FMI}\,($ \cal H$\,) is satisfied for some
 $ z \in \Bbb C \setminus \Bbb R $, then both \ 
$\mbox{\rm  TFMI}_{\mbox{\rm \scriptsize I}}\,(\cal H)$
and
$\mbox{\rm  TFMI}_{\mbox{\rm \scriptsize II}}\,(\cal H)$
are satisfied for the same $z$ as well.}
\end{lemma}

The opposite  implications (with respect to (\ref{7.3}), (\ref{7.4}))
 may be false, because the matrices $M_1(z)$ and $M_2(z)$
 are not invertible: $e_0\,T=0$, and the matrix $M_2(z)$ is not even  square.
Actually,
\begin{equation}
\mbox{FMI}\,(\cal H)\, \ \
 \mbox{$\not\Rightarrow $}\ \
\mbox{\rm TFMI}_{\mbox{\scriptsize\rm I}}\,(\cal H)
\label{7.5}
\end{equation}
Indeed, the product in the left hand side does not contain the
$nn$-th entry  $s_{2n}$ of the matrix $A$ at all, and 
 the positivity
of the matrix $A$
(and hence, the positivity of the matrix of the FMI\,($\cal H$)\,)
 depends essentially on this entry. However, the
FMI\,($\cal H$)\ {\sl and the
$\mbox{\rm TFMI}_{\mbox{\rm I}}$\,($\cal H$)\ are
``almost equivalent''
}: the matrix $M_{1}(z)$ (\ref{6.9}) is ``almost invertible''.
Since $T^{\ast}T=P$, where $P$ is a projector matrix:
$P=\mbox{\rm diag}\, [1,\,\ldots ,\,1,\, 0]$ \
($p_{kk}=1,\, k=0,\,1,\,\ldots ,\, n-1;\, p_{nn}=0 $), then
\begin{equation}
 \left[
\begin{array}{ccc}
T^{\ast}(I-\bar{z}T)^{-1}&\hspace*{-0.1cm} & 0\cr
 \vspace*{-0.1cm}&\vspace*{-0.1cm}\hspace*{-0.3cm} &\vspace*{-0.1cm} \cr
0&\hspace*{-0.1cm} &1\cr
\end{array}
\right]
\cdot
M_1(z)=
 \left[
\begin{array}{ccc}
P_{n-1}& \ & 0\cr
 \vspace*{-0.1cm}&\vspace*{-0.1cm}\hspace*{-0.3cm} &\vspace*{-0.1cm} \cr
0&\hspace*{-0.3cm} &1
\end{array}
\right]\cdot
\label{7.6}
\end{equation}
Hence,
{\sl
 the inequality, which is obtained from the inequality (\ref{6.1})
by replacing\footnote{
The last inequality is nothing more than the FMI of the form (\ref{6.1}), 
which is constructed from the ``truncated'' date 
$s_0,\, s_1,\, \ldots,\, s_{n-2}$. (The FMI (\ref{6.1}) is
constructed from the data $s_0,\, s_1,\, \ldots,\, s_{2n}$.)
}
 the matrix$A$ by the matrix $PAP$  and the column $B_{w}(z)$ by the column
$PB_{w}(z)$, holds.
}

{\sl
The inequalities FMI\,($\cal H$) and 
$
\mbox{\rm TFMI}_{\mbox{\scriptsize\rm II}}\,(\cal H)
$
are equivalent,
}
because there exists a left inverse matrix to the
matrix $M_{}(z)$: 
\begin{equation}
N(z)=
\left[
\begin{array}{ccc}
I &\hspace{-0.2cm} & 0  \cr
\vspace*{-0.1cm}      & \vspace*{-0.1cm} & \vspace*{-0.1cm}  \cr
0 &\hspace{-0.2cm} & e_0\,(I-\bar{z}T)\cr
\end{array}
\right], \qquad
N(z)\cdot M_{2}(z)=
\left[
\begin{array}{ccc}
I &\hspace{-0.2cm} & 0  \cr
\vspace*{-0.1cm}      & \vspace*{-0.1cm} & \vspace*{-0.1cm}  \cr
0 &\hspace{-0.2cm} & 1\cr
\end{array}
\right].
\label{7.7}
\end{equation}
Thus, we have  proved that
\begin{equation}
\mbox{FMI}\,(\cal H)\, \ \ \Leftrightarrow \ \ 
\mbox{\rm TFMI}_{\mbox{\scriptsize\rm II}}\,(\cal H).
\label{7.8}
\end{equation}
More formally:
\vspace*{5pt}

\begin{lemma}
{\sl
The inequality {\rm FMI\,($\cal H$)} is satisfied at some point $z \in
\Bbb C \setminus \Bbb R$ if and only if the inequality
$
\mbox{\rm TFMI}_{\mbox{\scriptsize\rm II}}\,(\cal H)
$
is satisfied for the same $z$.
}
\end{lemma}

The matrix of the
$
\mbox{\rm TFMI}_{\mbox{\scriptsize\rm I}}\,(\cal H)
$
is invariant with respect to the change
$z \rightarrow \bar{z}$. Thus:

{\sl If the inequality 
$
\mbox{\rm TFMI}_{\mbox{\scriptsize\rm I}}\,(\cal H)
$
is satisfied at some point
 $z \in \Bbb C \setminus \Bbb R$,  than it is satisfied also at the conjugate
point $\bar{z}$}.

 The following statement is not so evident:
\vspace*{5pt}

\begin{lemma}
{\sl
If the {\rm FMI\,($\cal H$)} is satisfied at some point
 $z\in \Bbb C \setminus \Bbb R$,  than it is satisfied also at the conjugate
point $\bar{z}$.
}
\end{lemma}

PROOF. The {\rm FMI\,($\cal H$)} can be written in the form 
$$
\left[
\begin{array}{ccc}
(I-zT)A(I-\bar{z}T^{\ast}&\hspace{0.1cm} & u\cdot w(z)-v\cr
\vspace*{-3pt}        & \vspace*{-3pt}      & \vspace*{-3pt}\cr
w^{\ast}(z)\cdot u^{\ast}-v^{\ast}& \hspace{0.1cm}&
\displaystyle
\frac{w(z)-w^{\ast}(z)}{z-\bar{z}}
\end{array}
\right]
\geq 0.
$$
The claim of the lemma follows from the matrix identity
\begin{eqnarray}
\left[
\begin{array}{cc}
I&  (\bar{z}-z)u  \cr
 \vspace*{-4pt}  &\vspace*{-4pt}  \cr
0  &   1
\end{array}
\right]
\left[
\begin{array}{ccc}
(I-zT)\,A\,(I-\bar{z}T^{\ast})&\hspace{3pt}&u\cdot w - v \cr
\vspace*{-4pt} &\vspace*{-4pt} & \cr
w^{\ast}\cdot u^{\ast}- v^{\ast}&\hspace{3pt}&
\displaystyle
\frac{w(z)-w^{\ast}(z)}{z-\bar{z}}
\end{array}
\right]
\left[
\begin{array}{cc}
I&  0\cr
 \vspace*{-4pt}  &  \vspace*{-4pt} \cr
(z-\bar{z})u &   1
\end{array}
\right]\nonumber  \\
\ \nonumber  \\
=
\left[
\begin{array}{cc}
(I-\bar{z}T)\,A\,(I-zT^{\ast})&  u\cdot w^{\ast} - v \cr
\vspace*{-4pt} &\vspace*{-4pt} \cr
w\cdot u^{\ast}- v^{\ast}&
\displaystyle
\frac{w(z)-w^{\ast}(z)}{z-\bar{z}}
\end{array}
\right].\cr
\hspace{5pt}
\label{7.9}
\end{eqnarray}
(where $w$ is an arbitrary complex number; we have to put $w=w(z)$, then 
$w^{\ast}=w(\bar{z})$).
To obtain the identity (\ref{7.9}), we perform the matrix multiplication
and use the identity
\begin{equation}
(I-\bar{z}T)\,A\,(I-zT^{\ast})-(z-\bar{z})\,(u\cdot v^{\ast}
-v\cdot u^{\ast})=(I-zT)\, A \,(I-\bar{z}T^{\ast}),
\label{7.10}
\end{equation}
which is equivalent to the Fundamental Identity (\ref{1.11}).
\hfill $\boxx$\\[0.2cm]
Now we turn to the extraction of  interpolation information from the
FMI\,($\cal H$). 
\vspace*{0.1cm}

\begin{theo} ({\cmss From the FMI\,($\cal H$)  to the moment conditions}).
{\sl
Let the interpolation data
 $s_0,\,s_1 ,\,\ldots\, , \,s_{2n-1},\,s_{2n}$ for the
Hamburger moment problem be given. Let $w$ be a function of the class
\Rgot\,($\Bbb H$)
 and  let the
{\rm FMI\ ($\cal H$) (\ref{6.1})} for this $w$ be satisfied at every point
$z$ in the upper half plane.  Then the function $w$ is 
representable in the form 
$w=w_{\sigma}$ for some (uniquely determined) measure $\sigma$. This
measure satisfies the moment conditions (\ref{4.5}); the interpolation 
conditions (\ref{4.19}) are satisfied as well.
}
\end{theo}

PROOF. According to Lemma 7.3, the {\rm FMI\ ($\cal H$) is satisfied for every 
$z\in \Bbb C \setminus \Bbb R$. By Lemma 7.2, the
\mbox{${\rm TFMI}_{\mbox{\scriptsize\rm II}}$}\,($\cal H$)
 is satisfied for every
$z\in \Bbb C \setminus \Bbb R$. First, from the
\mbox{${\rm TFMI}_{\mbox{\scriptsize\rm II}}$}\,($\cal H$)
we obtain the positivity condition
\begin{equation}
\frac{W_{w}(z)-W^{\ast}_{w}(z)}
{z-\bar{z}}\geq 0 \qquad (\forall z \in \Bbb C \setminus \Bbb R).
\label{7.11}
\end{equation}
Secondly, we derive the estimate
\begin{equation}
y\,W_{w}(iy)=O(1)  \qquad (\mbox{\sl as}\ \  y \uparrow \infty ).
\label{7.12}
\end{equation}
According to the matrix version of Nevanlinna's theorem,
the matrix function $W_{w}(z)$ is representable in the form
\begin{equation}
W_{w}(z)=\int\limits_{\Bbb R}\frac{d\Sigma (\lambda)}{\lambda-z}
\qquad (\forall z \in \Bbb C \setminus \Bbb R),
\label{7.13}
\end{equation}
where $d\Sigma (\lambda)$ is a nonnegative matrix-valued measure and the
integral
\begin{equation}
s_0\,({\Sigma})=\int\limits_{R}d\Sigma (\lambda)
\label{7.14}
\end{equation}
exists in the proper sense. Moreover,
\begin{equation}
\lim_{y\ \uparrow \ \infty}-iy\, W_{w}(iy)= s_0\,({\Sigma}).
\label{7.15}
\end{equation}
From the 
\mbox{${\rm TFMI}_{\mbox{\scriptsize\rm II}}$}\,($\cal H$) (\ref{7.2}) 
(for $z=iy,\, y \to \infty$) and from (\ref{7.15}) it now follows, that
\begin{equation}
A - s_0\,({\Sigma}) \geq 0.
\label{7.16}
\end{equation}
Of course,the condition (\ref{1.4.a}) for $w$ (see (\ref{6.15})) follows from
 the inequality (\ref{6.15}). Thus, $w=w_{\sigma}$ for
some $\sigma:\, s_{0}\,(\sigma)< \infty$. Let us clarify the structure of
the measure $d\Sigma$\,. We can expect that 
$W_{w}=W_{\sigma}$, and hence (see (\ref{6.6})) that
\begin{equation}
d\Sigma(\lambda)=(I-\lambda T)^{-1}u
\cdot
 d\sigma (\lambda)
\cdot
u^{\ast}(I-\lambda T^{\ast})^{-1}.
\label{7.17}
\end{equation}
This is the case indeed. To prove (\ref{7.17}), we turn to the formula
(\ref{6.11}). The functions 
\mbox{$(I-zT)^{-1}$ and  $(I-zT^{\ast})^{-1}$} are
holomorphic near the real axis (actually, these function are entire).
Applying {\cmss the generalized Sieltjes inversion formula}
 (\cite{KaKr}, \S 2) 
to (\ref{6.11}), we obtain (\ref{7.17}). In particular 
(see (\ref{5.6}) and (\ref{7.17})), the equality
\begin{equation}
s_{0}\,(\Sigma)=A_{\sigma}
\label{7.18}
\end{equation}
holds. Now (\ref{7.16}) takes the form
\begin{equation}
A-A_{\sigma}\geq 0.
\label{7.19}
\end{equation}
 The inequality (\ref{7.19}) itself ensures the 
condition (\ref{4.5}.ii), but it  does not ensure the condition
 (\ref{4.5}.i).
However, we can also exploit the asymptotics
(\ref{7.15}). Taking into account the
concrete structure (\ref{6.13}) of the matrix-function $W_{w}$, we see
that the asymptotic (\ref{7.15}) together with (\ref{4.5}.ii) 
leads to the condition (\ref{4.19}).  From (\ref{4.19}) of course follow
the moment condition (\ref{4.5}.i).

Another way to obtain these results is to multiply the
equality (\ref{6.11}) by $(I-zT)$ from the left and by $(I-zT^{\ast})$
from the right and then upon comparing the asymptotics of  both sides,
we see that
\begin{equation}
T\,(A-A_{\sigma})T^{\ast}=0.
\label{7.20}
\end{equation}
Thus, the nonnegative matrix $A-A_{\sigma}$ vanishes at all vectors from the
image of the matrix $T$. The orthogonal complement to this image is
generated by the $(n+1)\times 1$ vector
\begin{equation}
e_n=
\left[
\begin{array}{ccccc}
0&0&\cdots &0 &1\cr
\end{array}
\right].
\label{7.21}
\end{equation}
Hence,
\begin{equation}
A=A_{\sigma}+\rho \cdot e_n^{\ast}e_n, \quad 
\mbox{where $\rho$ is a nonnegative number}.
\label{7.22}
\end{equation}
In view of (\ref{5.1}) and (\ref{5.7}), the representation (\ref{7.22})
is equivalent to the moment conditions (\ref{4.5}).
\hfill $\boxx$

\begin{rem}
To obtain the estimate for the function $b_{w,2n}$,
we could restrict ourself to the subinequality of the inequality
(\ref{7.2}):
\begin{equation}
\left[
\begin{array}{ccc}
s_{2n}& & b_{w,2n}(z)\cr
\vspace*{-0.1cm}&\vspace*{-0.1cm} &\vspace*{-0.1cm} \cr
 b^{\ast}_{w,2n}(z)& &\displaystyle
\frac{b_{w,2n}(z)- b^{\ast}_{w,2n}(z)}{z-\bar{z}}
\end{array}
\right]
\geq 0.
\end{equation}
We can obtain this inequality from the inequality(\ref{7.2}), by ``framing''
it with the matrix
$$
\left[
\begin{array}{ccccccccc}
0 & \cdots &0 &1 & \vdots & 0 & \cdots & 0&0\cr
\vspace{-0.1cm}&\vspace{-0.1cm} &\vspace{-0.1cm} &\vspace{-0.1cm} &\vspace{-0.1cm}  &
\vspace{11pt} &\vspace{-0.1cm} &\vspace{-0.1cm} &\vspace{-0.1cm} \cr
0 & \cdots &0 &0 & \vdots & 0 & \cdots & 0&1\cr
\end{array}
\right]\cdot
$$
Combining this with   (\ref{6.22}), we obtain the following ``truncated''
transformation:
\begin{equation}
m(z)\cdot
\left[
\begin{array}{ccc}
A& &B_{w}(z)\cr
\vspace*{-0.1cm} &\vspace*{-0.1cm} & \vspace*{-0.1cm}\cr
B^{\ast}_{w}(z)& &\displaystyle\frac{w(z)-w^{\ast}(z)}{z-\bar{z}}\cr
\end{array}
\right]
\cdot m^{\ast}(z)=
\left[
\begin{array}{ccc}
s_{2n}& & b_{w,2n}(z)\cr
\vspace*{-0.1cm}&\vspace*{-0.1cm} &\vspace*{-0.1cm} \cr
 b^{\ast}_{w,2n}(z)& &\displaystyle
\frac{b_{w,2n}(z)- b^{\ast}_{w,2n}(z)}{z-\bar{z}}
\end{array}
\right]\, ,
\label{7.24}
\end{equation}
where
\begin{equation}
m(z)=
\left[
\begin{array}{ccccccc}
0 & 0& \cdots &0  & 1 & \vdots  & 0 \cr
 \vspace{-0.0cm} &\vspace{-0.1cm} & \vspace{-0.1cm}
  &\vspace{-0.1cm} &\vspace{-0.1cm} &\vspace{-0.1cm} & \vspace{-0.1cm}\cr
\bar{z}^{n-1} &\bar{z}^{n-2}   &  \cdots & 1 &  0 &\vdots  &\bar{z}^{n} \cr
\end{array}
\right]\cdot
\label{7.25}
\end{equation}
A transformation of the FMI of approximately the form (\ref{7.25})
appeared in the paper \cite{Kov} by I.Kovalishina (see pages 460-461 of the
Russian original or pages 424-425 of the  English translation).
(I.Kovalishina used a step by step algorithm, and did not
introduce the matrix  (\ref{7.25}) explicitly, but it is possible to
extract this matrix from her considerations.) Starting from\footnote
{The paper \cite{Kov} was published in 1983 only, but author was aware
of its content
 much earlier.}
\cite{Kov}, the author considered  transformations of the
FMI for various problems on integral representations,both
  discrete and continuous in \cite{K2}. 
The nontruncated transformation 
$
\mbox{FMI}\,(\cal H)\, \rightarrow  
\mbox{\rm TFMI}_{\mbox{\scriptsize\rm II}}\,(\cal H)
$
was considered by author in \cite{K3}. Such a transformation was considered
also by T.Ivanchenko and L.Sakhnovich \cite{IS1}, \cite{IS2}.
The nontruncated transformation 
$
\mbox{FMI}\,(\cal H)\, \rightarrow  
\mbox{\rm TFMI}_{\mbox{\scriptsize\rm I}}\,(\cal H)
$
was considered (for other classes of functions  and in different
notation) in \cite{KKY}.
 Systematic development of 
transformations of the FMI was  also presented in the preprint \cite{K4},
but \cite{K4} is not easily available.
\end{rem}
 
\vspace*{0.2cm}
\begin{minipage}{15.0cm}
\section{\hspace{-0.4cm}.\hspace{0.2cm}TRANSFORMATION OF FMI
 \boldmath$(\, \cal NP\,)$\unboldmath .}
\end{minipage}\\[-0.2cm]
\setcounter{equation}{0}

It is very easy to extract  interpolation information from the
FMI\,($\cal NP$). For this goal we need not  transform the FMI.
However, we have already learnd that such transformations and related 
structures are objects which are interesting in themselves. Therefore,
we will discuss transformations of the FMI\,($\cal NP$).
 (We know, that to a large extent 
 such transformations  depend only on Fundamental Identity for 
the considered problem and  not on the concret expression for the
entries in this
identity.) Thus, we consider a FMI of the form (\ref{1.5}) with $B_{w}$
and  $C_{w}$ of
the forms (\ref{1.8}), (\ref{2.2}), (\ref{2.3})
 and (\ref{1.6}), respectively, and we assume, that the 
Fundamental Identity (\ref{1.10}) is satisfied.

Let the function $w$ which appears in FMI\,($\cal NP$) be of  the form
$w=w_{\sigma, c}$ as in (\ref{3.1}). To guess formulas for transformations
of the FMI, we first consider the matrix function
\begin{equation}
W_{\sigma}(z)
=\int\limits_{\Bbb T}(tI-T)^{-1}
\cdot
\frac{1}{2}\,\frac{t+z}{t-z}\,d\sigma(t)
\cdot
(\bar{t}I-T^{\ast})^{-1},
\label{8.1}
\end{equation}
which is obtained by inserting the Schwarz kernel into the formula (\ref{3.4}) for $A_{\sigma}$.
We would like to obtain $W_{\sigma}$ from $A_{\sigma}$. For this goal we 
use the identity
\begin{equation}
\frac{1}{2}\,\frac{T+zI}{T-zI}\,(tI-T)^{-1}=
\frac{1}{2}\,\frac{t+z}{t-z}\,(tI-T)^{-1}+\frac{z}{z-t}\,(zI-T)^{-1}, 
\label{8.2}
\end{equation}
which was constructed with formulas (\ref{3.4}) and (\ref{3.6})
for $A_{\sigma}$ and $B_{\sigma}$ in mind. Now we multiply the identity
(\ref{8.2}) by 
$
u\cdot d\sigma(t)\cdot u^{\ast}(\bar{t}I-T^{\ast})^{-1}
$
and integrate over $\Bbb T$. Taking into account (\ref{3.4}) and (\ref{3.6}),
we obtain
\begin{equation}
W_{\sigma}(z)=\frac{1}{2}\,\frac{T+zI}{T-zI}\,A_{\sigma}-
(zI-T)^{-1}u\cdot B^{\ast}_{\sigma,c}(1/\bar{z}).
\label{8.3}
\end{equation}
The last formula is a heuristic reason for the following

\begin{defi}
{\sl
Given a Hermitian matrix $A$, a matrix $T$ and vectors $u$ and $v$ which
satisfy the FI {\rm (\ref{1.10})}, we associate with each
 function  $w$, which is holomorphic in
 $\Bbb C \setminus \Bbb T$ and satisfies the symmetry condition
{\rm (\ref{1.1})},
the function $W_{w}$:
\begin{equation}
W_{w}(z)=\frac{1}{2}\,\frac{T+zI}{T-zI}\,A-
(zI-T)^{-1}u\cdot B^{\ast}_{w,c}(1/\bar{z}).
\label{8.4}
\end{equation} 
or, in detail,
\begin{eqnarray}
 W_{w}(z)=\frac{1}{2}\,\frac{T+zI}{T-zI}\,A +
 (zI-T)^{-1}u\cdot v^{\ast}\,(z^{-1}\,I-T^{\ast})^{-1}\cr
     \vspace{-0.3cm}                \cr
+(zI-T)^{-1}u\cdot w(z) \cdot u^{\ast}\,(z^{-1}\,I-T^{\ast})^{-1}.
\label{8.5}
\end{eqnarray}
}
\end{defi}

Using the FI (\ref{1.10}), we obtain also another representation for
 $W_{w}(z)$:
\begin{equation}
 W_{w}(z)=\frac{1}{2}\,A\,\frac{I+zT}{I-zT}+B_{w}(z)\cdot 
u\,\frac{z}{I-zT^{\ast}},
\label{8.6}
\end{equation}
or, in detail,
\begin{eqnarray}
 W_{w}(z)=\frac{1}{2}\,A\,\frac{I+zT}{I-zT} -
 (zI-T)^{-1}v\cdot u^{\ast}\,(z^{-1}\,I-T^{\ast})^{-1}\cr
     \vspace{-0.3cm}                \cr
+(zI-T)^{-1}u\cdot w(z) \cdot u^{\ast}\,(z^{-1}\,I-T^{\ast})^{-1}.
\label{8.7}
\end{eqnarray}
In other words:
\vspace*{5pt}

\begin{lemma}
{\sl
The matrix-function $W_{w}$ satisfies the symmetry condition
\begin{equation}
W_{w}(z)=-W^{\ast}_{w}(1/\bar{z}) 
\qquad (\forall z \in \Bbb C \setminus \Bbb T).
\label{8.8}
\end{equation}
}
\end{lemma}
\vspace*{-0.3cm}}

Using the FI (\ref{1.10}), we obtain also the following result:

\begin{lemma}
{\sl
The matrix-function $W_{w}$ satisfies the  identity
\begin{equation}
W_{w}(z)-TW_{w}(z)T^{\ast}=
u\cdot \varphi^{\ast}_{w}\,(1/\bar{z})-\varphi_{w}\, (z)\cdot u^{\ast},
\label{8.9}
\end{equation}
\vspace*{-0.1cm}
where
\vspace*{-0.1cm}
\begin{equation}
\varphi_{w}\,(z)=\frac{1}{2}\,\frac{T+zI}{T-zI}\,(u\cdot w(z)-v).
\label{8.10}
\end{equation}
}
\end{lemma}

\begin{rem}{\sl
For $z=0$, the expression on the left hand side of (\ref{8.9}) is equal to 
\mbox{$\frac{1}{2}\, (A-TAT^{\ast})$,} and the expression on the right  hand
 side is equal to
\mbox{$\frac{1}{2}\,(u\cdot v^{\ast}+v\cdot u^{\ast})$.}
 Thus, the formula (\ref{8.9}) is in some sense 
an analytic continuation of the FI (\ref{1.10})}
\end{rem}

\begin{rem}{\sl
The equality (\ref{8.9}), considered as an equation with respect to the
matrix $W_{w}(z)$, can be used to calculate this matrix.}
\end{rem}

Let us calculate the matrix $W_{w}(z)$ for the $\cal NP$ problem
with  data given by (\ref{2.2}) and (\ref{2.3}).
From the equation (\ref{8.9}),
we obtain the following formula:
\begin{equation}
W_{w}(z)=
\frac{1}{2}
\left\|
\begin{array}{c}
\underline{\displaystyle 
\frac{z_k+z}{z_k-z}\,(w_k-w(z))+
\frac{1+z\bar{z}_l}{1-z\bar{z}_l}\,(w(z)+w^{\ast}_l)} \cr
\displaystyle 1- z_k \bar{z}_l\cr
\end{array}
\right\|_{1\leq k,l\leq n}
\cdot
\label{8.11}
\end{equation} 
Let us introduce the matrices
\begin{equation}
M_{1}(z)=
\left[
\begin{array}{ccc}
(I-\bar{z}T)^{-1}&\ & \bar{z}\,(I-\bar{z}T)^{-1}u \cr
\vspace{-0.1cm}& \vspace{-0.1cm}& \vspace{-0.1cm}\cr      
0                   &\ &   1 \cr
\end{array}
\right]
\label{8.12}
\end{equation}
and
\begin{equation}
M_{2}(z)=
\left[
\begin{array}{ccc}
I                   &\ &   0 \cr
\vspace{-0.1cm}& \vspace{-0.1cm}& \vspace{-0.1cm}\cr      
\,(I-\bar{z}T)^{-1}&\ &\bar{z}\,(I-\bar{z}T)^{-1}u\cr
\end{array}
\right]\cdot
\label{8.13}
\end{equation} 
Performing the matrix multiplication, we obtain (after some calculations with
the  entries):
\begin{equation}
M_{1}(z)\cdot
\left[
\begin{array}{ccc}
A &B_{w}(z)\cr
& \cr
B^{\ast}_{w}(z)&\displaystyle\frac{w(z)-w^{\ast}(z)}{z-\bar{z}}\cr
\end{array}
\right]
\cdot
M^{\ast}_{1}(z)=
\left[
\begin{array}{ccc}
\displaystyle\frac{W_{w}(z)+W^{\ast}_{w}(z)}{1-z\bar{z}}& 
  &\displaystyle\frac{B_{w}(z)-B_{w}(1/\bar{z})}{1-z\bar{z}}\cr
&  & \cr
\displaystyle\frac{B^{\ast}_{w}(z)-B^{\ast}_{w}(1/\bar{z})}{1-z\bar{z}}&
 &\displaystyle\frac{w(z)+w^{\ast}(z)}{1-z\bar{z}} \cr
\end{array}
\right]
\label{8.14}
\end{equation}
and
\begin{equation}
M_{2}(z)\cdot
\left[
\begin{array}{ccc}
A &B_{w}(z)\cr
& \cr
B^{\ast}_{w}(z)&\displaystyle\frac{w(z)-w^{\ast}(z)}{z-\bar{z}}\cr
\end{array}
\right]
\cdot
M^{\ast}_{2}(z)=
\left[
\begin{array}{ccc}
\displaystyle A& 
  &\displaystyle W_{w}(z)+\frac{A}{2}\cr
       & & \cr
\displaystyle W^{\ast}_{w}(z) +\frac{A}{2} & &
\displaystyle \frac{W_{w}(z)+W^{\ast}_{w}(z)}{1-z\bar{z}}\cr
\end{array}
\right]\cdot
\label{8.15}
\end{equation}
The calculations mentioned above  are based essentially
on the following consequence of the FI (\ref{1.10}):
\vspace*{0.2cm}
\begin{eqnarray}
(z-T)^{-1}\,A\,(\bar{z}-T^{\ast})^{-1} \hspace{10.0cm}\cr
\  \ \cr
=\frac{1}{1-z\bar{z}}\,
\left(
\frac{1}{2}\,\frac{T+zI}{T-zI}\,A+
\frac{1}{2}\,A\,\frac{T^{\ast}+\bar{z}I}{T^{\ast}-\bar{z}I}
\right)+
(zI-T)^{-1}\cdot
\frac{u\, v^{\ast}+v\,u^{\ast}}{1-z\bar{z}}
\cdot (\bar{z}I-T^{\ast})^{-1}. \cr
\label{8.16}
\end{eqnarray}
We consider two variants of theTransformed Fundamental Matrix Inequality
 (for the
Nevanlinna-Pick problem):
the $\mbox{\rm TFMI}_{\mbox{\scriptsize\rm I}}(\cal NP)$ and
the $\mbox{\rm TFMI }_{\mbox{\scriptsize\rm II}}(\cal NP)$.
 
$\mbox{\rm TFMI}_{\mbox{\scriptsize\rm I}}(\cal NP)$ has the form
\begin{equation}
\left[
\begin{array}{ccc}
\displaystyle\frac{W_{w}(z)+W^{\ast}_{w}(z)}{1-z\bar{z}}& 
  &\displaystyle\frac{B_{w}(z)-B_{w}(1/\bar{z})}{1-z\bar{z}}\cr
&  & \cr
\displaystyle\frac{B^{\ast}_{w}(z)-B^{\ast}_{w}(1/\bar{z})}{1-z\bar{z}}&
 &\displaystyle\frac{w(z)+w^{\ast}(z)}{1-z\bar{z}} \cr
\end{array}
\right]
\geq 0.
\label{8.17}
\end{equation}
$\mbox{\rm TFMI}_{\mbox{\scriptsize\rm II}}(\cal NP)$ has the form
\begin{equation}
\left[
\begin{array}{ccc}
\displaystyle A& 
  &\displaystyle W_{w}(z)+\frac{A}{2}\cr
       & & \cr
\displaystyle W^{\ast}_{w}(z) +\frac{A}{2} & &
\displaystyle \frac{W_{w}(z)+W^{\ast}_{w}(z)}{1-z\bar{z}}\cr
\end{array}
\right]
\geq 0.
\label{8.18}
\end{equation}
We see that both of theseTFMI's contain the function $W_{w}(z)$.

\begin{defi}
{\sl
Given a $\cal NP$ problem with interpolation nodes
 $z_1 ,\, z_2,\, \ldots ,\,z_n$
 in the unit disc $\Bbb D$,
 the point $z \in \Bbb C \setminus \Bbb T$ is said to be
{\cmss nonsingular} ,
if $z \not= 0,\,\infty;\, z_1,\,z_2,\,\ldots,\,z_n;\,
\bar{z}^{-1}_1,\,\bar{z}^{-1}_2,\,\ldots,\,\bar{z}^{-1}_n$.
}
\end{defi}

If $z$ is a nonsingular point, then the matrices
 $(zI-T)^{-1},\,(I-\bar{z}T)^{-1}$
are defined (and, of course, invertible). (Strictly speaking, we can define
the  matrices
$W_{w}(z)$, $M_{1}(z)$ and $M_{2}(z)$ for nonsingular $z$ only).
 {\sl  For nonsingular $z$,
the matrix $M_{1}(z)$ is invertible and the matrix $M_{2}(z)$ has
 a left inverse.}

\begin{lemma}
{\sl
Let $z \in \Bbb C \setminus \Bbb T$ be a nonsingular point. Then
the FMI\,($\cal NP$) is
satisfied at this point $z$ if and only if each of two inequalities
$\mbox{\rm TFMI}_{\mbox{\scriptsize\rm I}}(\cal H)$
and
$\mbox{\rm TFMI}_{\mbox{\scriptsize\rm II}}(\cal H)$
is satisfied at this point.
}
\end{lemma}

\begin{lemma}
{\sl
Let $z \in \Bbb C \setminus \Bbb T$ be a nonsingular point. Then the
 FMI\,($\cal NP$) is
satisfied at this point $z$ if and only if it is satisfied at the ``symmetric'' point
$\bar{z}^{-1}$ as well.
}
\end{lemma}

PROOF. The FMI\,($\cal NP$) is equivalent to the inequality
\begin{equation}
\left[
\begin{array}{ccc}
(zI-T)\,A(\bar{z}I-T^{\ast})& 
  &u\cdot w(z)-v\cr
       & & \cr
w^{\ast}(z)\cdot u^{\ast}-v^{\ast} & &
\displaystyle \frac{w_{w}(z)+w^{\ast}_{w}(z)}{1-z\bar{z}}\cr
\end{array}
\right]
\geq 0.
\label{8.19}
\end{equation}
The claim of the lemma follows from the matrix identity
\begin{eqnarray}
\left[
\begin{array}{cc}
I &-(1-z\bar{z})\,u  \cr
     &  \cr
0    & 1
\end{array}
\right]\hspace{-2pt}
\left[
\begin{array}{cc}
(I-zT)\,A\,(I-\bar{z}T^{\ast})& u\cdot w - v \cr
\vspace*{-2pt} &\vspace*{-2pt} \cr
w^{\ast}\cdot u^{\ast}- v^{\ast}&
\displaystyle
\frac{w(z)+w^{\ast}(z)}{1-z\bar{z}}
\end{array}
\right]\hspace{-2pt}
\left[
\begin{array}{cc}
I & 0\cr
   &    \cr
-(1-z\bar{z})\,u^{\ast}     &   1
\end{array}
\right]=  \nonumber \\
\ \nonumber \\
\left[
\begin{array}{ccc}
(I-\bar{z}T)\,A\,(I-zT^{\ast})& & u\cdot w^{\ast} - v \cr
\vspace*{-2pt} &\vspace*{-2pt} & \vspace*{-2pt}\cr
w\cdot u^{\ast}- v^{\ast}&\hspace{3pt}&
\displaystyle
\frac{w(z)+w^{\ast}(z)}{1-z\bar{z}}
\end{array}
\right], \nonumber \\
\label{8.20}
\end{eqnarray}
where  $w=w(z)$ and $w^{\ast}=-w(1/\bar{z})$.
 To obtain the identity (\ref{8.20}), we perform 
the matrix multiplication and use the identity
\begin{equation} 
(zI-T)\,A\,(\bar{z}I-T^{\ast})-
(1-z\bar{z})\,(u\cdot v^{\ast}+v\cdot u^{\ast})
=(I-\bar{z}T)\, A \,(I-zT^{\ast}),
\label{8.21}
\end{equation}
which is equivalent to the FI (\ref{1.10}).	

\begin{lemma}
{\sl
The $\mbox{\rm TFMI}_{\mbox{\scriptsize\rm II}}(\cal NP)$ (\ref{8.18})
 holds for every point
 $z \in \Bbb D$ if and only if the function $W_{w}(z)$
satisfies the positivity condition:
\begin{equation}
W_{w}(z)+W^{\ast}_{w}(z) \geq 0  \qquad (\forall z \in \Bbb D).
\label{8.22}
\end{equation}
}
\end{lemma}
\vspace*{-2pt}
PROOF. The implication 
$
\mbox{\rm TFMI}_{\mbox{\scriptsize\rm II}} \Rightarrow \mbox{\rm(\ref{8.22})}
$
is evident. The opposite implication is nothing more that the Schwarz-Pick 
inequality for the function $W_{w}(z)$ for the points: $0$ and $z$
(because $W_{w}(0)=\frac{A}{2}$).
\hfill $\boxx$

From Lemmas 8.3 and 8.5 we obtain the following conclusion:

\begin{theo}
{\sl
A function $w$, holomorphic in $\Bbb C \setminus \Bbb T$
 and satisfying the symmetry condition (\ref{1.1}),
 satisfies the FMI\,($\cal NP$)
for all $z \in \Bbb D$ (or, what the same for all 
$z \in \Bbb C \setminus \Bbb T$) 
if and only if the function $W_{w}(z)$ which is defined by (\ref{8.4})
 satisfies the positivity
condition (\ref{8.22}).
}
\end{theo}

Taking into account the concrete form (\ref{8.11})
 of the matrix $W_{}$ for the $\cal NP$ problem, we obtain:

\begin{theo}
{\sl
Let  the interpolation data for the $\cal NP$ 
problem (\ref{2.1}) in the function class 
\Cgot\,($\Bbb D$) be given by  (\ref{2.2}) and (\ref{2.3}). 
A function $w$, which is holomorphic in $\Bbb D$, is a solution of
the $\cal NP$ problem (with these data) if and only if the real part of the
matrix on the right hand side of (\ref{8.11}) is nonnegative for every
$z \in \Bbb D$. 
}
\end{theo}

\begin{rem}
{\sl
The matrix in (\ref{8.11}) is an orthogonal projection of the
operator \mbox{$\displaystyle \frac{1}{2}(I+zU)\,(I-zU)^{-1}$,} where $U$ is a
generalised unitary extension of some isometric  operator,
related to the considered problem. 
}
\end{rem}

This is a consequence of the
 $\mbox{\rm TFMI}_{\mbox{\scriptsize\rm II}}(\cal NP)$. A consecuence
of the $\mbox{\rm TFMI}_{\mbox{\scriptsize\rm I}}(\cal NP)$
also may be interesting. The inequality (\ref{8.17}) is equivalent to the
 inequality
\begin{equation}
\left[
\begin{array}{cc}
W_{w}(z) + W^{\ast}_{w}(z)&B_{w}(z)-B_{w}(1/\bar{z}) \cr
\vspace*{-4pt} &\vspace*{-4pt} \cr
B^{\ast}_{w}(z)-B^{\ast}_{w}(1/\bar{z})& w(z) + w^{\ast}(z)\cr
\end{array}
\right]
\geq 0
\qquad (\forall z \in \Bbb D)\, .
\label{8.23}
\end{equation}
The matrix function on the left hand side of (\ref{8.23}) is harmonic and 
nonnegative in $\Bbb D$  and hence it admits a Riesz-Herglotz
representation. Let
\begin{equation}
\left[
\begin{array}{cc}
d\Sigma (t)& d\mu (t)\cr
\vspace*{-4pt} &\vspace*{-4pt} \cr
d\mu^{\ast} (t)& d\sigma(t)\cr
\end{array}
\right]
\label{8.24}
\end{equation}
be the block decomposition of the representing measure. Now we can apply
\u{S}mul'yan's results
 from\footnote{The paper \cite{S} by
Yu.L.\,\u{S}mul'yan looks as it was written especially to be used in
this paper.} \cite{S}, to obtain the inequality
\begin{equation}
\int\limits_{\Bbb T}\,d\mu(t)\,(d\sigma(t))^{-1}\,d\mu^{\ast}(t)\leq
\int\limits_{\Bbb T}\,d\Sigma(t)\,,
\label{8.25}
\end{equation}
where the integral on the left hand side is the so called 
{\sl Operator Hellinger Integral}. Because
\vspace{-4pt}
$$
W_{w}(0)+W^{\ast}_{w}(0))=A,\quad \mbox{\rm it follows that}
\int\limits_{\Bbb T}\,d\Sigma = A.
\vspace{-4pt}
$$
Thus
\begin{equation}
\int\limits_{\Bbb T}\,d\mu(t)\,(d\sigma(t))^{-1}\,d\mu^{\ast}(t)\leq A\,.
\label{8.26}
\end{equation}
It is not difficult to show that in the considered case (the $\cal NP$ 
problem with finitely many interpolation nodes located inside $\Bbb D$)
the equality holds in (\ref{8.26}). 
In the general situation, $A$ is a nonnegative Hermitian form in some
 vector space.
Then, the $\mbox{\rm TFMI}_{\mbox{\scriptsize\rm I}}(\cal NP)$
leads to the representation of a nonnegative Hermitian form 
by the {\cmss   Hellinger Integral}. It is worthy to mention
that it was the Hellinger  integral, which was
used for the integral representation of Hermitian kernels early
in the development of the theory. In more recent time, the Stieltjes integral
ousted  the Hellinger integral from this circle of problem.
 However, the use of the Stieltjes integral
leads to difficulties. It may not exist, and we have to use
rigged Hilbert spaces and all that. And the Hellinger integral exists
always (and under some conditions it may be reduced to the Stieltjes
integral). By our opinion, the use of the Hellinger integral lies in the
essence of matter. The moral is clear:

{\cmss GO BACK TO THE CLASSICS.}
\newpage
$$$$
$$$$

\vspace{0.6cm}
\begin{minipage}{8.2cm}
Victor Katsnelson \\
Department of Theoretical Mathematics\\
The Weizmann Institute of Science\\
Rehovot, 76100\\ Israel\\
e-mail: katze@wisdom.weizmann.ac.il
\end{minipage}
\vskip .5cm      
{ AMS subject classification:} 30D50, 46E10.

\end{document}